\documentclass[12pt]{amsart}
\usepackage{lscape}

\usepackage[ansinew]{inputenc}
\usepackage{amsmath,euscript,amssymb,amsthm,amscd,graphicx,epsfig}
\usepackage{amsmath,euscript,amssymb,amsthm,xypic,amscd,graphicx,epsfig}

\usepackage[all]{xy}

\newtheorem{thm}{Theorem}[section]
\newtheorem{cor}[thm]{Corollary}
\newtheorem{proposition}[thm]{Proposition}
\newtheorem{lem}[thm]{Lemma}

\newtheorem{definition}[thm]{Definition}

\title{Homology of spaces of regular loops in the sphere}

\author{David Chataur and Jean-Fran\c cois Le Borgne}

\address{Laboratoire Paul Painlev\'e,
Universit\'e de Lille1, 59655 Villeneuve d'Ascq C\'edex, France}

\email{David.Chataur@math.univ-lille1.fr}

\email{Leborgne.Jean-Francois@math.univ-lille1.fr}

\begin{document}

\maketitle

\begin{abstract}
In this paper we compute the singular homology of the space of immersions
of the circle into the $n$-sphere. Equipped with Chas-Sullivan's loop
product these homology groups are graded commutative algebras, we also compute these algebras. 
We enrich Morse spectral sequences
for fibrations of free loop spaces together with loop products, this
offers some new computational tools for string topology. 
\end{abstract}

\section*{Introduction}

The aim of this paper is to compute the singular homology with integral
coefficients of the spaces $Imm(S^{1},S^{n})$ of immersions of the
circle $S^{1}$ into the $n$-sphere $S^{n}$ also called regular loops. These spaces play a
key role in knot theory because they detect non-trivial homology classes
for spaces of embeddings $Emb(S^{1},S^{n})$. For example in \cite{CCL} the authors
use a desingularization map in order to produce non-trivial cohomology
classes for spaces of knots $Emb(S^{1},\mathbb{R}^{n})$ from cohomology
classes of spaces of singular knots i.e. spaces of immersions with
a fixed number of transverse double points.

Let us explain the strategy of the computation of these homology groups. We consider the space
of immersions $Imm(S^{1},S^{n})$ as an open subset of the free loop
space $\mathcal{L}S^{n}=C^{\infty}(S^{1},S^{n})$ and we endow the
$n$-sphere with its standard Riemannian metric. We consider the energy
functional $E$ defined by
\[
\begin{array}{ccc}
\mathcal{L}S^{n} & \rightarrow & \mathbb{R}\\
\gamma & \mapsto & \int_{S^{1}}\left\Vert \dot{\gamma}(\theta)\right\Vert ^{2}d\theta\end{array}\]
and we also consider its restriction $E_{imm}$ to the space of immersions. From this functional we get
a filtration of the space of immersions by the energy level. It follows
from Morse theory that this filtration is well understood in the case of the
space $\mathcal{L}S^{n}$ of all smooth loops and gives a very efficient way to compute its singular homology. 
But it is a priori more difficult to handle it for
immersions. To this filtration one naturally associates a spectral sequence
wich converges to $H_{*}(Imm(S^{1},S^{n}),\mathbb{Z})$. In order
to compute the $E_{1}$-term of this spectral sequence it is essential
to understand the homology of the pairs $(Imm^{<\lambda_{p}+\epsilon},Imm^{<\lambda_{p}-\epsilon})$
where $\lambda_{p}$ is a critical value of $E_{imm}$ and \[
Imm^{<a}:=\left\{ \gamma\in Imm(S^{1},S^{n})/E_{imm}(\gamma)<a\right\} .\]

Let $US^{n}$ be the unit tangent bundle of $S^{n}$ and $D$ be
the map \[
\begin{array}{ccc}
Imm(S^{1},S^{n}) & \rightarrow & \mathcal{L}US^{n}\\
\gamma & \mapsto & \frac{\dot{\gamma}}{\left\Vert \dot{\gamma}\right\Vert }\end{array}\]
by a famous result of Hirsch and Smale \cite{Sm} this map is a weak homotopy
equivalence. We have the following factorization \[
E_{imm}=E\circ\pi\circ D:Imm(S^{1},S^{n})\stackrel{D}{\longrightarrow}\mathcal{L}US^{n}\stackrel{\pi}{\longrightarrow}\mathcal{L}S^{n}\stackrel{E}{\longrightarrow}\mathbb{R}\]
where $\pi=\mathcal Lp$ is the loop map associated to the canonical projection $p:US^{n}\rightarrow S^{n}$.
Thus our problem reduces to a problem of ''fiberwise'' Morse theory,
we will consider the filtration by the energy level on the space $\mathcal{L}US^{n}$
induced by the map $E\circ\pi$. From this filtration one gets a Morse-Serre
type spectral sequence. To be more precise one has a fibration \[
\mathcal{L}S^{n-1}\rightarrow\mathcal{L}US^{n}\rightarrow\mathcal{L}S^{n}\]
and the filtration on the total space is induced by the Morse filtration
(the energy filtration) of the base. We will show that this fiberwise filtration is homotopy equivalent to the filtration induced by $E_{imm}$.

If we consider this spectral sequence just as a spectral sequence of abelian groups it seems impossible to complete the
computation, we need to enrich this spectral sequence with an additional
algebraic data. Here string topology enters the game. From the foundational
work of Chas and Sullivan one knows that 
$$\mathbb H_*(\mathcal LUS^n):= H_{*+2n-1}(\mathcal{L}US^{n})$$
is a graded commutative algebra equipped with the so-called loop product \cite{CS}.
Hingston and Goresky have proved that this
spectral sequence is multiplicative for the loop product \cite{HG}. 

Moreover we prove that this spectral sequence collapses at the $E_{2}$-term ($E_2=E_{\infty}$). Thus
we are left with some non trivial extensions issues that we solve by comparing the preceding spectral sequences with the
Serre spectral sequences associated to the fibrations
\[
\Omega US^{n}\rightarrow\mathcal{L}US^{n}\rightarrow US^{n}\]
 \[
\mathcal{L}S^{n-1}\rightarrow\mathcal{L}US^{n}\rightarrow\mathcal{L}S^{n}.\]
In fact the Serre spectral sequence associated to the first fibration collapses at the $E_2$-term when $n$ is even while the second collapses at the $E_2$-term when $n$ is odd (\cite{L}). But in each cases we encounter extensions issues that we are able to solve by comparison with the Morse spectral sequence. Thus if we want to complete the computation we need to use these three spectral sequences together. Then we get our main result :

\begin{thm} \label{main}
Let $n \geq 2$ be an integer.
Then $\mathbb{H}_*(Imm(S^1, S^{2n})) \simeq \mathbb{H}_*(\mathcal LUS^{2n})$ is isomorphic to the algebra
$$\mathbb{Z}[y_{-2n}, \alpha_{2n-2}, \beta_{4n-2}]\otimes \Lambda(x_{-4n+1},k_{-1})/\mathcal{R}$$
where $\mathcal{R}$ is the ideal generated by
$$( 2y_{-2n},2k_{-1}, 2\alpha_{2n-2},x_{-4n+1}y_{-2n}, x_{-4n+1}k_{-1}, y_{-2n}^2, y_{-2n}k_{-1}-x_{-4n+1}\alpha_{2n-2})$$
and $\mathbb{H}_*(Imm(S^1, S^{2n+1})) \simeq \mathbb{H}_*(\mathcal LUS^{2n+1})$  
is isomorphic to the algebra 
$$\mathbb{Z}[v_{2n}, y_{-2n}, u_{4n-2}] \otimes \Lambda(x_{-2n-1},\theta_{-1})/(y^2_{-2n},\theta_{-1} y_{-2n}, 2u_{4n-2}y_{-2n}).$$
\end{thm}

{\bf Remarks} Let us give some additional informations about the preceding computations. 
\\
(1) First let us notice that the evaluation map
$$ev_0:Imm(S^1,S^n)\rightarrow US^n$$
has a section $geod$ this section associates to each unit tangent vector $u\in US^n$ the unique great circle (the prime geodesic) $\gamma_u$ such that $\dot{\gamma_u}(0)=u$. In homology $(ev_0)_*$ and $(geod)_*$ are morphisms of algebra between the intersection homology algebras $\mathbb H_*(US^n)$ and $\mathbb H_*(Imm(S^1,S^n))$. Explicitely one can identify the algebra 
$\mathbb H_*(US^{2n})$ with the subalgebra of $\mathbb H_*(Imm(S^1,S^{2n}))$ generated by the classes $y_{-2n}$ and $x_{-4n+1}$, in the odd case we can identify $\mathbb H_*(US^{2n+1})$ with  the subalgebra of $\mathbb H_*(Imm(S^1,S^{2n+1}))$ generated by the classes $y_{-2n}$ and $x_{-2n-1}$.
\\
\\
(2) Let $Imm_{u}(S^1,S^n)=ev_0^{-1}(u)$ be a fiber of the evaluation map 
$$ev_0:Imm(S^1,S^n)\rightarrow US^n.$$ 
We consider the inclusion
$$j:Imm_{u}(S^1,S^n)\rightarrow Imm(S^1,S^n)$$
in homology we have the intersection morphism
$$I=j_!:\mathbb H_*(Imm(S^1,S^n))\rightarrow H_*(Imm_{u}(S^1,S^n))$$ 
which corresponds to the classical intersection morphism
$$Int:\mathbb H_*(\mathcal LUS^n)\rightarrow H_*(\Omega US^n)$$
of string topology. In the even case, we have
$$H_*(\Omega US^{2n})\cong \mathbb Z[\alpha_{2n-2},\beta_{4n-2}]/(2\alpha_{2n-2})$$
the morphism $I$ is given by $I(x_{-4n+1})=I(y_{-2n})=I(k_{-1})=0$
and $I(\alpha_{2n-2})=\alpha_{2n-2},I(\beta_{4n-2})=\beta_{4n-2}$. In the odd case we have
$$H_*(\Omega US^{2n+1})\cong H_*(\Omega S^{2n+1})\otimes H_*(\Omega S^{2n})\cong \mathbb Z[v_{2_n}]\otimes \mathbb Z[f_{2n-1}]$$
and $I(x_{-2n-1})=I(y_{-2n})=I(\theta_{-1})=0$, $I(v_{2n})=v_{2n}$, $I(u_{4n-2})=(f_{2n-1})^2$.
\\
\\
(3) Over the field of rational numbers these computations become easier because rationaly $US^{2n}\simeq_{\mathbb Q} S^{4n-1}$ and
$US^{2n+1}\simeq_{\mathbb Q} S^{2n+1}\times S^{2n}$. 
\\
\\
(4) The computation of the loop algebra $\mathbb H_*(\mathcal LUS^n)$ is also related to two other topological problems. One knows that the space $US^n$ is homotopy equivalent to the configuration space $F_3(S^n)$ of three points into $S^n$, in fact  the projection
$$proj_{3,1}:F_3(S^n)\rightarrow S^n$$
given by $proj_{3,1}(x_1,x_2,x_3)=x_1$ is fiberwise homotopy equivalent to
$$p:US^n\rightarrow S^n.$$
The computation of the homology of $\mathcal LUS^n$ is essential into the understanding of probems of 3-body type into $S^n$, we refer the reader to Fadell and Husseini's monograph \cite{FH}. This computation is related to symplectic topology, from recent work of Abbondandolo and Schwarz one knows that the loop algebra $\mathbb H_*(\mathcal LUS^n)$ is isomorphic as an algebra to Floer homology $HF_*((TUS^n)^*)$ of the cotangent bundle of $US^n$ together with the pair of pants product \cite{AS}.
\\
\\
Plan of the paper :
\\
Section 1 : we recall some basic facts about immersion spaces and string topology.
\\
Section 2 : we build our main technical tool the Morse-Serre spectral sequence and show its compatibity with the loop product. We play with various filtrations of loop spaces : the length filtration, the energy filtration and the filtration by the square root of the energy.
\\
Section 3 : this section is devoted to the computation of the 0-th column of the Morse-Serre spectral sequence.  
\\
Section 4 : we compute $\mathbb H_*(Imm(S^1,S^{2n}))$.
\\
Section 5 : we compute $\mathbb H_*(Imm(S^1,S^{2n+1}))$.
\\
\\
{\bf Acknowledgement.} The authors thank Fran\c cois Laudenbach for a careful reading of a first version of this paper. 
The authors are supported by the ``Laboratoire Paul Painlev\'e, UMR 8524
de l'Universit\'e des Sciences et Technologie de Lille et du CNRS'', by the GDR
2875 ``Topologie Alg\'ebrique et Applications du CNRS''. The first author's research
is supported in part by ANR grant 06-JCJC-0042 ``Op\'erades, Big\`ebres et
Th\'eories d'Homotopie''.
 
\section{String topology of immersion spaces}

\subsection{Loop product}

We recall the definition of Chas-Sullivan's loop poduct \cite{CS}.
Let $M$ be a $d$-dimensional connected, compact oriented manifold. Moreover we suppose that $M$ is closed that is to say without boundary. 
Let $\delta_{M}:M\hookrightarrow M\times M$, $x\mapsto(x,x)$ be
the diagonal embedding. We also denote by $\widetilde{\delta}_{M}:\mathcal{L}M\times_{M}\mathcal{L}M\hookrightarrow\mathcal{L}M\times\mathcal{L}M$
the embedding of composable loops where\[
\mathcal{L}M\times_{M}\mathcal{L}M=\{(\gamma_{1},\gamma_{2})\in\mathcal{L}M\times\mathcal{L}M/\gamma_{1}(0)=\gamma_{2}(0)\}.\]
 We have the following pull-back diagram:

\[
\begin{array}{ccc}
\mathcal{L}M\times_{M}\mathcal{L}M & \stackrel{\widetilde{\delta}_{M}}{\longrightarrow} & \mathcal{L}M\times\mathcal{L}M\\
\downarrow &  & \qquad\quad\downarrow ev_{0}\times ev_{0}\\
M & \stackrel{\delta_{M}}{\longrightarrow} & M\times M\end{array}\]

and we have a composition map $comp_{M}:\mathcal{L}M\times_{M}\mathcal{L}M\to\mathcal{L}M$
that concatenates the composable loops. To be more precise, the composition map is a continuous
map given by $comp_M(\gamma,\gamma')(t)=\gamma(2t)$ if $0\leq t\leq 1/2$ and 
$comp_M(\gamma,\gamma')(t)=\gamma'(2t-1)$ if $1/2\leq t\leq 1$. 
\\
As a pull-back of the diagonal embedding, the embedding $\widetilde{\delta_{M}}$ is smooth
and finite codimensional, so it is possible to define a shriek map
\[
\widetilde{\delta_{M}}_{!}:H_{*}(\mathcal{L}M\times\mathcal{L}M)\to H_{*-d}(\mathcal{L}M\times_{M}\mathcal{L}M).\]

Let us brievely recall the construction of this shriek map for a smooth
embedding $j:X\hookrightarrow Y$ of finite codimension $k$, where
$X$ and $Y$ can be, and this is our case, infinite dimensional manifolds. We suppose that this embedding is co-oriented, that is to say the normal bundle of $j$ is an oriented vector bundle. For
such an embedding we suppose that there is a tubular neighbourhood $Tub(X)$ of $j(X)$
that is diffeomorphic to the total space of the disc bundle $D_k(X)\to X$
of the normal bundle of $j$ and let $S_{k-1}(X)\to X$ be its associated
sphere bundle. 
Then $j_{!}:H_{*}(Y)\to H_{*-k}(X)$ is defined as
the composition of:

a) The inclusion of pairs $H_{*}(Y)\to H_{*}(Y,Y-j(X))$

b) The excision isomorphism \[
H_{*}(Y,Y-j(X))\to H_{*}(Tub(X),Tub(X)-j(X))\]

c) An isomorphism induced by an homotopy equivalence \[
H_{*}(Tub(X),Tub(X)-j(X))\to H_{*}(D_{k}(X),S_{k-1}(X))\]

d) The Thom isomorphism $H_{*}(D_{k}(X),S_{k-1}(X))\to H_{*-k}(X)$
that is the composition of the cap product with the Thom class $\tau\in H_{k}(D_{k}(X),S_{k-1}(X))$
with the map $\pi_{*}$ induced by the canonical projection $\pi:D_{k}(X)\to X$. 
If the embedding is not co-oriented we have to use singular homology with local coefficients. 
\\
The Chas and Sullivan loop product $\mu$ is defined as the following
composition map:

\[
\mu:H_{*}(\mathcal{L}M)\otimes H_{*}(\mathcal{L}M)\stackrel{\times}{\to}H_{*}(\mathcal{L}M\times\mathcal{L}M)\stackrel{\widetilde{\delta_{M}}_{!}}{\to}H_{*-d}(\mathcal{L}M\times_{M}\mathcal{L}M)\stackrel{comp_{M_{*}}}{\to}H_{*-d}(\mathcal{L}M).\]

where $\times$ denotes the cross product. For other constructions of this product we refer the reader to Sullivan's survey paper \cite{S1}. Let us recall some basic facts about this product
\begin{enumerate}
    \item When suitably regraded i.e. we define $\mathbb H_*(\mathcal LM):=H_{*+d}(\mathcal LM)$ the loop product is unitary and commutative.
    \item The loop product is compatible with the intersection product of $M$. Let us recall that $\mathbb H_*(M):=H_{*+d}(M)$ together with the intersection product is a graded commutative algebra called the intersection algebra. The evaluation map $ev_0:\mathcal LM\rightarrow M$ and its section 
$c:M\rightarrow \mathcal LM$ the constant loop map induce two morphisms of algebras
$$(ev_0)_*:\mathbb H_*(\mathcal LM)\rightleftarrows \mathbb H_*(M):(c)_*.$$
The algebra $\mathbb H_*(M)$ is isomorphic to a sub-algebra of the loop algebra, the unit of $\mu$ is equal to $(c)_*([M])$ where $[M]$ is the fundamental class of $M$.
    \item The composition of based loops induces in homology
$$\sharp:H_*(\Omega M)\otimes H_*(\Omega M)\rightarrow H_*(\Omega M)$$
the Pontryagin product which has a unit and is associative. The inclusion
$j:\Omega M\rightarrow \mathcal LM$ can be considered as a codimension $d$ embedding. Thus one can define a shriek map
$$Int:=j_!:\mathbb H_*(\mathcal LM)\rightarrow H_*(\Omega M).$$ 
This is a morphism of algebras called the intersection morphism.
\item There are several constructions of the loop product, some are purely algebraic, others use stable homotopy theory. A chain level description of the loop product has been given by F. Laudenbach in \cite{La}, the author uses his approach to give a nice construction of the multiplicative structures considered in \cite{HG} and in this paper.   
\end{enumerate}

\subsection{Hirsch-Smale's theorem and loop products}
Let $M$ be a smooth Riemannian manifold and let $UM$ be the unit tangent bundle. We consider the evaluation map
$$ev_0:Imm(S^1,M)\rightarrow UM$$
defined by $ev_0(\gamma)=\frac{\dot{\gamma}(0)}{\|\dot{\gamma}(0)\|}$. This evaluation map is a fibration, this was first proved by S. Smale \cite[Theorem B]{Sm1}. But one can prove a stronger result, in fact this map is locally trivial. We use a theorem of R. S. Palais \cite[Theorem A]{Pal} : the group of diffeomorphisms $Diff(M)$ acts on $Imm(S^1,M)$ and $UM$, the evaluation map is equivariant and $UM$ is a $Diff(M)$-space admitting local cross-sections (let $x\in UM$ be a unitary tangent vector there is a map $\chi$ of a neighborhood U of $x$ into $Diff(M)$ such that $\chi(u)(x)=u$ for all $u\in U$), then by \cite[Theorem A]{Pal} the map $ev_0$ is locally trivial.
\\ 
 We have the following pull-back diagram:
\[
\begin{array}{ccc}
Imm(S^1,M)\times_{UM} Imm(S^1,M) & \stackrel{\widetilde{\delta}_{UM}}{\longrightarrow} & Imm(S^1,M)\times Imm(S^1,M)\\
\downarrow &  & \qquad\quad\downarrow ev_{0}\times ev_{0}\\
UM & \stackrel{\delta_{UM}}{\longrightarrow} & UM\times UM\end{array}\]
and a composition map 
$$comp:Imm(S^1,M)\times_{UM} Imm(S^1,M)\rightarrow Imm(S^1,M)$$ 
one can proceed as in the preceding section and define an immersion product and an immersion algebra
$$\mathbb H_*(Imm(S^1,M)):=H_{*+2d-1}(Imm(S^1,M)).$$
We also notice that if $u\in UM$ and if $Imm_u(S^1,M):=(ev_0)^{-1}(u)$ we have an intersection morphism
$$\mathbb H_*(Imm(S^1,M))\rightarrow H_*(Imm_u(S^1,M)).$$ 

{\bf Remark.} The definition of the map $comp$ needs to be modified slightly. First we notice that the concatenation of two composable immersions is not well defined. Let $\gamma$ and $\gamma'$ be two composable immersions. At $t=0$ the loops $\gamma$ and $\gamma'$ do not have the same tangent vector but their normalizations are the same by definition of the space of composable immersions. Thus after concatenation of the composable immersions we reparametrize the loop by its arc length. We also notice that the concatenation is a piecewise $C^2$-path. This space of piecewise $C^2$-immersions is homotopy equivalent to space $Imm(S^1,M)$ of $C^{\infty}$-immersions. We should also work with piecewise $C^2$-paths rather than $C^{\infty}$-paths (all these spaces of loops and paths are homotopy equivalent see \cite{Pal2}).
\\
\\
Now we relate the immersion algebra to the string topology intersection algebra of $UM$. By Hirsch-Smale's homotopy theory of immersions (\cite{H}, \cite{Sm}) we know that the differential 
$$D:Imm(S^1,M)\rightarrow \mathcal LUM$$ 
is a homotopy equivalence. We have a map of fibrations
\[
\begin{array}{ccccc}
Imm_u(S^1,M)&\stackrel{j}{\longrightarrow}&Imm(S^1,M) & \stackrel{ev_0}{\longrightarrow} & UM\\
\downarrow D_u&&\downarrow D&  & \downarrow Id\\
\Omega_u UM&\stackrel{j}{\longrightarrow}&\mathcal LUM & \stackrel{ev_0}{\longrightarrow} & UM.\end{array}\]  
This map is a morphism of fiberwise monoids in the sense of Gruher and Salvatore \cite{GS}. Thus we have:

\begin{lem}
Let $M$ be a connected closed smooth manifold then $D_*$ gives an isomorphism of algebras between the immersion algebra 
$\mathbb H_*(Imm(S^1,M))$ and the loop algebra $\mathbb H_*(\mathcal LUM)$. 
\end{lem}

Let us suppose that $M$ is $1$-connected. Following Cohen-Jones-Yan, the Serre spectral sequence associated to $ev_0$ is multiplicative with respect to the loop product. This structure involves the intersection product on the homology of the base and the Pontryagin product on the homology of the fiber. Our first computational tool is this multiplicative spectral sequence:
$$\mathbb E^2_{p,q}:=\mathbb H_p(UM,H_q(\Omega UM))\Rightarrow \mathbb H_{p+q}(\mathcal LUM).$$
Let us consider the unit tangent bundle 
$$S^{d-1}\rightarrow UM\stackrel{p}{\longrightarrow} M.$$
Our second computational tool is associated to the fibration $\pi:=\mathcal Lp$
$$\mathcal LS^{d-1}\rightarrow \mathcal LUM\stackrel{\pi}{\longrightarrow} \mathcal LM.$$
The second's author has proved in his PhD thesis \cite{L} that the Serre spectral sequence of this fibration is multiplicative for the loop product. We have
$$\mathbb E^2_{p,q}:=\mathbb H_p(\mathcal LM,\mathbb H_q(\mathcal LS^{d-1}))\Rightarrow \mathbb H_{p+q}(\mathcal LUM).$$
If $U_m\subset M$ is a neighbourhood of a point $m\in M$ diffeomorphic to $\mathbb R^d$, the differential maps the space $imm(S^1,U_m)$ into $\mathcal LU\mathbb R^{d}$. Let us denote by $triv$ the trivilization map 
$$triv:U\mathbb R^{d}\rightarrow \mathbb R^d\times S^{d-1}$$
and we define $tr:\mathcal LU\mathbb R^d\rightarrow \mathcal LS^{d-1}$ as $tr:=\mathcal L(pr_2\circ triv)$. If 
$f:Imm(S^1,M)\rightarrow \mathcal LUM$ is the canonical inclusion,
then we have a diagram
\[
\begin{array}{ccccc}
Imm(S^1,U_m)&\stackrel{}{\longrightarrow}&Imm(S^1,M) & \stackrel{f}{\longrightarrow} & \mathcal LM\\
\downarrow tr\circ D&&\downarrow D&  & \downarrow Id\\
\mathcal  LS^{d-1}&\stackrel{}{\longrightarrow}&\mathcal LUM & \stackrel{\pi}{\longrightarrow} & \mathcal LM.\end{array}\]  
In this diagram each vertical map is a weak homotopy equivalence. The right hand side square is strictly commutative whereas the left hand side square is commutative only up to homotopy.

\subsection{Loop algebra for spheres} We have seen that the loop algebra of spheres plays a major role in the determination of the immersion algebra of a manifold. Let us recall Cohen-Jones-Yan's computation of these algebras
\[
\mathbb{H}_{*}(\mathcal{L}S^{n})\cong\left\{ \begin{array}{c}
\Lambda(a)\otimes\mathbb{Z}[u]\quad for\: n\: odd\\
(\Lambda(b)\otimes\mathbb{Z}[a,v])/(a^{2},ab,2av)\quad for\: n\: even\end{array}\right.\]
In the Serre spectral sequence of
$$\Omega S^n\rightarrow \mathcal LS^n \rightarrow S^n$$
we have $a\in\mathbb H_{-n}(\mathcal LS^n)\cong \mathbb E^{\infty}_{-n,0}$, 
$b\in \mathbb H_{-1}(\mathcal LS^n)\cong \mathbb E^{\infty}_{-n,n-1}$,
$u\in \mathbb H_{n-1}(\mathcal LS^n)\cong \mathbb E^{\infty}_{0,n-1}$ and
$v\in \mathbb H_{2n-2}(\mathcal LS^n)\cong \mathbb E^{\infty}_{0,2n-2}$.
In the next section we will show how to recover these results using Morse theory.

\section{Loop Product and Morse Theory}

\subsection{Morse theory for loop spaces}

In this section we recall some basic facts about Morse-Bott theory for loop spaces, we refer the reader to Bott's papers \cite{B} and to \cite{Bnd}.
Let $W$ be a manifold and $f$ be a real valued function on $W$.
A connected submanifold
$N\subset W$ is said to be a nondegenerate critical manifold of $W$ if the
following conditions are satisfied.

1) Each point $p\in N$ is a critical point of $f$.

2) The Hessian of $f$ is non degenerate in the normal direction to
$N$.

Spelled out this last condition takes this form: 

3) we consider a small tubular $\epsilon$-neighborhood
$W_{\epsilon}(N)$ of $N$, which is fibered over $N$ by the normal
discs swept out by geodesics of lenght $\leq\epsilon$ in the normal
direction to $N$, relative to some Riemannian structure on $W$. Then,
the condition 3) is equivalent to the following assumption: 
\\
\textit{$f$ restricted
to each normal disc is nondegenerate.} 
\\
In this case we decompose the normal bundle $\nu N$ into a positive
and a negative part. \[
\nu N=\nu^{+}N\oplus\nu^{-}N,\]
 where $\nu_{p}^{+}N$ and $\nu_{p}^{-}N$ are respectively spanned
by the positive and negative eigen-directions of the Hessian of $f$.
The fiber dimension of $\nu^{-}N$ will be denoted $\alpha^{N}$ and
refered as the index of $N$ rel $f$.
We refer to this conditions as ``the Bott non degeneracy conditions'' and say that $f$ is Morse-Bott.
\\
Let us concentrate on the case of the free loop space $\mathcal LM$ of a $d$-dimensional compact Riemannian manifold $M$.
We work with piecewise smooth loops, we consider the Sobolev class $H^1(S^1,M)$. Such a model has the advantage to be a Hilbert manifold.
The function we consider is the energy functional
\[
E:\mathcal LM\to\mathbb{R},\gamma\mapsto\int_{S^{1}}{\Vert \dot{\gamma}(t)\Vert}^{2}dt.\]
The critical points of the energy are the closed geodesics. In order to use Morse theory we pick a metric on $M$ such that this function is Morse-Bott. Thus the critical points of $E$ are collected on compact
critical manifolds that satisfy the Bott non-degeneracy conditions. 
\\
In general this condition holds for a generic metric on $M$. According to R. Bott \cite{B} this genericity follows already from general position arguments in Morse's work.
In W. Ziller
\cite[Theorem 2]{Z} proves that globally symmetric spaces satisfy this condition. This result was proved precedently by R. Bott for the loop space with fixed endpoints $\Omega(M,p,q)$ of a globally symmetric space.
\\
Let $0=\lambda_{0}<\lambda_{1}<...$
be the critical values of the energy function $E$ and 
$$\mathcal LM^{\leq\lambda_{i}}:=E^{-1}(]-\infty;\lambda_{i}]).$$
The spaces $\mathcal LM^{\leq\lambda_{i}}$ give rise to a filtration of $\mathcal LM$ that
provides a filtration of the singular chain complex $C_{*}(\mathcal LM)$.
We define \[
F_{p}C_{p+q}(\mathcal LM):=C_{p+q}(\mathcal LM^{\leq\lambda_{p}}).\]
 Now let us identify the graded module  associated to this filtration. 
Let $\Sigma_{r}$ be the critical submanifold
associated to the critical value $\lambda_{r}$ namely $$\Sigma_{r}=\{\gamma\in \mathcal LM/E(\gamma)=\lambda_{r},\mbox{ and }dE(\gamma)=0\},$$ 
this manifold is finite dimensional. 
\\
For the moment we assume that the critical submanifold is connected.
The tangent bundle of $\mathcal LM$ restricted to $\Sigma_{r}$ splits in
three parts: \[
T\mathcal LM_{\vert_{\Sigma_{r}}}\simeq\mu^{-}_{r}\oplus\mu^{0}_{r}\oplus\mu^{+}_{r}\]
 corresponding to the signature of the hessian of $E$ at the point
$\gamma$. The main result of Morse theory for free loop spaces is that there is a homotopy equivalence
$$\mathcal LM^{\leq \lambda_r}\simeq\mathcal LM^{\leq \lambda_{r-1}}\cup_{f_r}D(\mu^{-}_{r})$$
for a gluing map $f_r:S(\mu^{-}_{r})\rightarrow \mathcal LM^{\leq \lambda_{r-1}}$. 
The quotient $\mathcal LM^{\leq\lambda_{r}}/\mathcal LM^{\leq\lambda_{r-1}}$
is homotopically equivalent to the Thom space $Th(\mu^{-}_{r})$
of the bundle $\mu^{-}_{r}$. This proves that 
$$F_{p}C_{p+*}(\mathcal LM)/F_{p-1}C_{p+*}(\mathcal LM)=C_{p+*}(\mathcal LM^{\leq\lambda_{p}},\mathcal LM^{\leq\lambda_{p-1}})$$
is quasi-isomorphic to $\widetilde{C}_{p+*}(Th(\mu^{-}_{p}))$. The case of 
several isolated critical manifolds can be treated in the same
fashion.
\\
The filtration leads to a homology spectral sequence called the Morse spectral sequence. 

\begin{thm}
\label{Morse spec seq} Let $M$ be a metric such that the energy is a Morse-Bott function. The energy filtration of $C_{*}(\mathcal LM)$ induces
a spectral sequence 
$$\{ E_{*,*}^{r}(\mathcal{M})(\mathcal LM)\}_{r\in\mathbb{N}}$$
converging to $H_{*}(\mathcal LM)$ 
\[
E_{p,q}^{r}(\mathcal{M})(\mathcal LM)\Rightarrow H_{p+q}(\mathcal LM).\]
We suppose that the critical set $\Sigma_p$ with critical value $\lambda_p$ is a union of connected non-degenerate manifolds $N_i$  
and that $\mu^{-}_i$ denotes the negative part of $T\mathcal LM_{\vert_{N_{i}}}$.
Then the $E^{1}$-page $E_{p,q}^{1}(\mathcal{M})(\mathcal LM)=H_{p+q}(\mathcal LM^{\leq\lambda_{p}},\mathcal LM^{\leq\lambda_{p-1}})$
is isomorphic to the reduced
homology $\oplus_i\widetilde{H}_{p+q}(Th(\mu^{-}_i))$. 
\end{thm}

{\bf Remark.} Let $N_i$ be a connected component of a critical set $\Sigma_p$, the theorem above and the Thom isomorphism gives an isomorphism
$$H_*(\mathcal LM^{\leq\lambda_{p}},\mathcal LM^{\leq\lambda_{p-1}})\cong \widetilde{H}_{*-a_p}(N_i,\mathbb Z)$$
only when the negative bundle $\mu^{-}_i$ of dimension $a_p$ is oriented otherwise one has to use coefficients in the orientation bundle of $\mu^{-}_i$.

\subsection{Loop products in the Morse spectral sequence}
In this section we lift the Chas-Sullivan product at the chain level in order to get a multiplicative Morse spectral sequence. In \cite{La} F. Laudenbach has given an alternative construction of the multiplicative structure of this spectral sequence.  
\\
\\
{\bf The use of the length filtration} The composition map $comp_M$ is not compatible with the energy filtration, in \cite{HG} the authors changed this map by a parameterized composition map. However, if one wants to use the map $comp_M$ one has to work with the length filtration rather than the energy filtration. Because of its geometric flavour this filtration is more natural and will be very useful in the fiberwise case.  
\\
Let $\mathcal LM^{\leq a}_l$ be the space of loops $\gamma$ of length $L(\gamma)\leq a$. By the Cauchy-Schwarz inequality one has the inclusion 
$$\mathcal LM^{\leq a}\rightarrow \mathcal LM^{\leq a}_l$$ 
this inclusion is a homotopy equivalence and this map has a homotopy inverse. Let us be more precise we consider the 
 subspace $\mathcal A(M)^{\leq a}\subset \mathcal LM$ of loops $\gamma$ parametrized proportionnaly to arc length such that $L(\gamma)\leq a$. If $\gamma\in \mathcal A(M)^{\leq a}$ the length $L(\gamma)$ is equal to $F(\gamma)$. 
One can also filter the space $\mathcal A(M)$ by the length of loops by \cite[prop 2.2]{HG} this filtration is homotopy equivalent to the Morse filtration of $\mathcal LM$ by $F$ and to the length filtration of $\mathcal LM$. In fact the inclusion  $\mathcal A(M)\subset \mathcal LM$ has a homotopy inverse 
$$A:\mathcal LM\rightarrow \mathcal A(M)$$
which associates to any path the same
path parametrized proportionally to arclength, with the same basepoint. This map does not change the lentgh of the loop. We have the homotopy equivalences of pairs
$$(\mathcal LM^{\leq a},\mathcal LM^{\leq b})\leftarrow(\mathcal A(M)^{\leq a},\mathcal A(M)^{\leq b})\rightarrow(\mathcal LM_l^{\leq a},\mathcal LM_l^{\leq b}).$$
From now on, we use the length filtration in order to take advantage of the formula
$$L(comp_M(\gamma,\gamma'))=L(\gamma)+L(\gamma').$$
By abuse of notation we use the same notation for the length filtration and the energy filtration.
\\
\\
{\bf The chain level product} Let $j:X\rightarrow Y$ be a co-oriented finite codimensional embedding and let $C_{*}(X)$ be the singular chain complex of $X$. If we chose $\hat{\tau}\in C^{k}(D_{k}(X),S_{k-1}(X))$ representing the Thom
class at the chain level, we can define $\hat{j_{!}}:C_{*}(Y)\to C_{*-k}(X)$
that induces $j_{!}$ in homology whatever the choice of the cocycle
representing $\tau$ is.
So let us consider the Chas and Sullivan loop product at the chain level, we have
the following diagram
\[
\begin{array}{c}
C_{p+q}(\mathcal LM^{\leq\lambda_{p}})\otimes C_{p'+q'}(\mathcal LM^{\leq\lambda_{p'}})\\
\downarrow\times\\
C_{p+q+p'+q'}(\mathcal LM^{\leq\lambda_{p}}\times \mathcal LM^{\leq\lambda_{p'}})\\
\downarrow\widetilde{\delta_{M}}\\
C_{p+q+p'+q'-d}(\mathcal LM^{\leq\lambda_{p}}\times_{M}\mathcal LM^{\leq\lambda_{p'}})\\
\;\downarrow {(comp_M)_*}\\
C_{p+q+p'+q'-d}(\mathcal LM^{\leq\lambda_{p}+\lambda_{p'}})\end{array}\] 
\\
\\
{\bf Condition (Cl).} Let $M$ be as in the preceding section a d-dimensional compact oriented Riemannian manifold. The main hypothesis that we need here is that all geodesics $\gamma$ are closed and simply periodic with the same prime length $l$. We recall that this means that $\gamma(0) = \gamma(1)$, $\dot{\gamma}(0) = \dot{\gamma}(1)$, $\gamma$ is injective on $(0,1)$ and the length of $\gamma$ namely $L(\gamma)$ is equal to $l$ if $\gamma$ is prime. 
\\
Moreover we suppose that all the negative bundles of critical sets are oriented. We denote this condition by the condition (Cl). Spheres, complex projective spaces, more generally $1$-connected globally symmetric spaces of rank one all satisfy this condition. In fact R. Bott proved in \cite{BCl} that the singular cohomology of a space that satisfies condition (Cl) is isomorphic to the cohomology of an irreductible symmetric space of rank $1$ . 
\\
The fundamental point for us is that the condition (Cl) implies that $\lambda_r + \lambda_{r'} = rl + r'l = (r+r')l = \lambda_{r+r'}$.
Under this condition the loop product gives a product 
$$F_pC_{p+q}(\mathcal LM) \otimes F_{p'}C_{p'+q'}(\mathcal LM)\rightarrow F_{p+p'}C_{p+q+p'+q'-d}(\mathcal LM).$$ 
This implies the following proposition.
\begin{proposition} \label{mult}
Let $M$ be a $d$-dimensional compact Riemannian manifold that satisfies condition (Cl).
Then, the loop product induces a multiplicative structure on the shifted Morse spectral sequence $\mathbb{E}^r_{p,q}(\mathcal{M})(\mathcal LM) := E^r_{p, q+d}(\mathcal{M})(\mathcal LM)$.
\end{proposition}

{\bf Remark.} Of course the inequality $\lambda_r + \lambda_{r'} \leq \lambda_{r+r'}$ is essential in order to have a multiplicative structure. It would be interesting to investigate the existence of metrics that satisfy this condition, in particular in the case of connected sums of symmetric spaces.  
\\
\\
The following Lemma gives the bi-graded module structure of $\mathbb{E}^1_{*,*}(\mathcal{M})(\mathcal LM)$ in function of $H_*(M)$, $H_*(UM)$ and of $\alpha_r$. Where $\alpha_r$ is the index of the critical value $\lambda_r$ namely the dimension of the fiber of $\mu^{-}_r \to \Sigma_r$. 

\begin{lem} \label{lemma 1}
If $M$ is a Riemannian manifold that satisfies condition (Cl), then for $r>0$ we have isomorphisms of modules
$$\phi_r: H_*(\mathcal LM^{\leq \lambda_r}, \mathcal LM^{\leq \lambda_{r-1}}) \to H_{*- \alpha_r}(UM).$$
\end{lem}

Proof.
By retraction along the gradient flow lines of the energy functional the space $\mathcal LM^{< \lambda_{r}}$ retracts to  $\mathcal LM^{\leq \lambda_{r-1}}$. This gives an isomorphism  $$H_*(\mathcal LM^{\leq \lambda_r}, \mathcal LM^{\leq \lambda_{r-1}}) \to H_*(\mathcal LM^{\leq \lambda_r}, \mathcal LM^{< \lambda_{r}}).$$
The right term is equal to $H_*(\mathcal LM^{\leq \lambda_r}, \mathcal LM^{\leq \lambda_{r}}-\Sigma_r)$.
By excision isomorphism, it is isomorphic to $H_*(D_{\alpha_r}(\Sigma_r), S_{\alpha_r-1}(\Sigma_r))$.
Next, $H_*(D_{\alpha_r}(\Sigma_r), S_{\alpha_r-1}(\Sigma_r))$ is isomorphic to $H_{*- \alpha_r}(UM)$ by the Thom isomorphism.
Since $M$ satisfies condition (Cl), then $\Sigma_r$ is diffeomorphic to $UM$. The composition of these isomorphisms gives the required isomorphism $\phi_r: H_*(\mathcal LM^{\leq \lambda_r}, \mathcal LM^{\leq \lambda_{r-1}}) \to H_{*- \alpha_r}(UM).$
\rightline{$\square$}

\begin{cor} \label{cor 1}
 The bigraded module $\oplus_{p,q \geq 0} E^1_{p,q}(\mathcal M)(\mathcal LM)$ is isomorphic to $H_q(M)$ if $p=0$ and to $\oplus_{p\geq 1,q \geq 0} H_{p+q-\alpha_p}(UM)$ if $p \geq 1$.
\end{cor}
Proof. This is direct application of Theorem \ref{Morse spec seq} and of Lemma \ref{lemma 1}.
\rightline{$\square$}
\\
\\
{\bf The multiplicative structure.} In what follows we reformulate Corollary 12.7 of Theorem 12.5 of \cite{HG}. We introduce a bigraded algebra $\mathcal A_{*,*}$. As bigraded module, $\mathcal A_{*,*}$ is a regraduation by a translation of $d$ of the bigraded module of Lemma \ref{lemma 1}. We set $\mathcal A_{0,q}:= \mathbb H_q(M)$ and 
$\mathcal A_{p,q}:=\mathbb H_{p+q-\alpha_p}(UM)<T^p>$ if $p>0$ where $T$ is an element of bidegree $(1,\alpha_1+d-2)$. Moreover we suppose that we have Bott's iteration formula:
$$\alpha_p=p\alpha_1+(p-1)+(d-1).$$
The multiplicative structure of $\mathcal A_{*,*}$ is given by : 
\\
1) the intersection algebra $\mathbb H_{*}(M)$ if $p=0$, 
\\
2) the algebra $\mathbb H_{p + q -\alpha_{p}}(UM)[T]_{\geq 1}$ of polynomials of degree $\geq 1$ in $T$ when $p>0$, 
\\
3) the products involving an element of $\mathbb H_{*}(M)$ and an element of $\mathbb H_{p + q-\alpha_{p}}(UM)[T]_{\geq 1}$ is given by a topological $\mathbb H_*(M)$-module structure on $\mathbb H_*(UM)$ we denote it by 
$$m: \mathbb{H}_s(M) \otimes \mathbb{H}_t(UM) \to \mathbb{H}_{s+t}(UM).$$ 
This structure is given by the map
$$m= \hat{\delta_M}! \circ \times$$
where $\times$ is the cross product 
$$H_i(M)\otimes H_j(UM)\rightarrow H_{i+j}(M\times UM)$$
and $\hat{\delta_M}!$ is the gysin map associated to $\hat{\delta_M}$ defined using the pull-back map
$$
\xymatrix{
M \times_M UM \ar@{^{(}->}[r]^{\hat{\delta_M}} \ar[d] & M \times UM \ar[d]^{id \times p} \\
M \ar@{^{(}->}[r]^{\delta_M} & M \times M }
$$
at the homology level we have the morphism 
$$\mathbb H_u(M) \otimes \mathbb H_v(UM)T^k \to \mathbb H_{u+v}(UM)T^k$$
$$x \otimes yT^k \mapsto m(x,y)T^k.$$

\begin{thm} \label{morseE1}
Let $M$ be a $d$-dimensional compact Riemannian manifold that satisfies condition (Cl).
As algebra, $\mathbb{E}^1_{*,*}(\mathcal{M})(\mathcal LM)$ is isomorphic to $\mathcal A_{*,*}$.
In other words, the multiplicative structure on the shifted Morse spectral sequence $\mathbb{E}^r_{p,q}(\mathcal{M})(\mathcal LM) := E^r_{p, q+d}(\mathcal{M})(\mathcal LM)$ is given at the $E^1$-level by:
$$\mathbb{E}^1_{*,*}(\mathcal{M})(\mathcal LM)= \mathbb{H}_*(M) \oplus \mathbb{H}_*(UM)[T]_{\geq 1}.$$
The bidegree of $T$ is $(1, \alpha_1+d -2)$, elements of $\mathbb H_*(M)$ and of $\mathbb H_*(UM)$ are of bidegree $(0,*)$. We have 
$$\mathbb{E}^1_{0,q}(\mathcal{M})(\mathcal LM)= \mathbb{H}_q(M)$$
and for $p \geq 1$,
$$\mathbb{E}^1_{p,q}(\mathcal{M})(\mathcal LM)= \mathbb{H}_{q-p \alpha_1}(UM)<T^p>$$
The multiplication between the $0$-th column and the others is induced by the $\mathbb{H}_*(M)$-module structure on $\mathbb{H}_*(UM)$.
\end{thm}

Proof.
The Corollary \ref{cor 1} tells us that $\mathbb{E}^1_{*,*}(\mathcal{M})(\mathcal LM)$ and $\mathcal A_{*,*}$ are isomorphic as modules. 
It remains to prove that there multiplicative structures are isomorphic. This is a direct consequence of the diagram page 39 of \cite{HG} see also \cite{La}. The diagonal map on $UM$ factorizes as 
$$UM\stackrel{\delta_S}{\rightarrow} UM\times_M UM \stackrel{\delta_M}{\rightarrow} UM\times UM$$
For simplicity, we denote $\mathcal LM^{\leq \lambda_r}$ by $\Lambda_r$. At the level of homology we have the following commutative diagram

$$
\xymatrix{
H_p(\Lambda_r, \Lambda_{r-1}) \otimes H_{p'}(\Lambda_{r'}, \Lambda_{r'-1}) \ar[r]^{\phi_r \otimes \phi_{r'}} \ar[d]^{\times} & H_{p-\alpha_r}(UM) \otimes H_{p'-\alpha_{r'}}(UM)  \ar[d]^{\times} \\
H_{p+p'}(\Lambda_r \times \Lambda_{r'}, \Lambda_r \times \Lambda_{r'-1} \bigcup \Lambda_{r-1} \times \Lambda_{r'}) \ar[r]^{\phi_r \times \phi_{r'}} \ar[d]^{\widetilde{\delta_M}!} & H_{p-\alpha_r + p'-\alpha_{r'}}(UM \times UM)  \ar[d]^{\hat{\delta_{M}}!} \\
H_{p+p'-d}(\Lambda_r \times_M \Lambda_{r'}, \Lambda_r \times_M \Lambda_{r'-1} \bigcup \Lambda_{r-1} \times_M \Lambda_{r'}) \ar[r]^{\phi_r \times_M \phi_{r'}} \ar[d]^{{comp_M}_*} & H_{p-\alpha_r + p'-\alpha_{r'} -d}(UM \times_M UM)  \ar[d]^{\bar{\delta_S}!} \\
H_{p+p'-d}(\Lambda_{r +r'}, \Lambda_{r + r'-1}) \ar[r]^{\phi_{r+r'}} & H_{p-\alpha_r + p'-\alpha_{r'} -d -(d-1)}(UM)  }
$$
We recall Bott's iteration formula \cite{BCl}
$$\alpha_{r+r'} = \alpha_r + \alpha_{r'} + d-1.$$
Then we have: 
$$
\xymatrix{
H_{p+p'-d}(\Lambda_{r + r'}, \Lambda_{r + r'-1}) \ar[r]^{\phi_{r+r'}} & H_{p +p' -d -\alpha_{r + r'}}(UM).}
$$ 
\\
Thus we have for $p \geq 1$, the following multiplicative identification
$$\mathbb{E}^1_{p,q}(\mathcal{M})(\mathcal LM)= \mathbb{H}_{q-p \alpha_1}(UM)<T^p>.$$
For $p=0$, we observe that $\mathbb{H}_*(\mathcal LM^{\leq 0},\mathcal LM^{\leq -1})$ with the loop product is in fact $\mathbb{H}_*(M)$ with the intersection product. This allows us to identify the $0^{th}$ column $$\mathbb{E}^1_{0,q}(\mathcal{M})(\mathcal LM)= \mathbb{H}_q(M).$$
And for $p$ or $p'$ equal to zero, we get the $\mathbb{E}^1_{0, *}(\mathcal{M})$-module structure given by $m$ on $\mathbb{E}^1_{*, *}(\mathcal{M})(\mathcal LM)$.

\rightline{$\square$}

\subsection{The example of spheres}
Let us illustrate the preceding theorem, we consider the case of spheres. 
The computation of the Betti numbers was done by Bott in \cite{Bnd}. 
Bott's method was used by Ziller to completely determine the homology of free loop spaces of symmetric spaces of rank one \cite{Z}. Ziller proved that the Morse spectral sequence collapses at the $E_1$-term. 
\\
We recall that on the $n$-sphere together with its standard Riemannian structure, the closed geodesics are represented by non-degenerate critical manifolds of two types:
\\
$Geod_0$= $S^n$ corresponding to constant loops,
\\
$Geod_k$= $US^n$ corresponding to the great circles starting to a unit vector $u\in US^n$ and traversing that circle $k$ times.
\\
The indices for $k>0$ are given by the formula 
$$index(Geod_k)=(2k-1)(n-1).$$
\\
We consider the Morse spectral sequence for spheres. The $E^1$-term is generated by elements on $E^1_{0,*}$ and on $E^1_{1,*}$. As we have a spectral of algebra the differential $d^1$ is a derivation, thus to compute $d^1$ it suffices to compute it on $E^1_{0,*}$ and $E^1_{1,*}$. Thanks to inclusion of the constant loops in $\mathcal LS^n$, which is a section of the evaluation map, we deduce that all the elements of $E^1_{0,*}$ must survive. We thus have that $d^1=0$. As algebra we find that $E^1=E^2$ and for degree reasons the derivation satisfy $d^r=0$ for $r>1$. Finally we have $E^{\infty}=E^1$. 
\\
In the even case, there is no extension issue. In the odd case, one is left with one multiplicative extension issue, but one can solve it by using the intersection morphism 
$$\mathbb H_*(LS^n)\rightarrow H_*(\Omega S^n).$$  
\\
At the end we recover that 
\[
\mathbb{H}_{*}(\mathcal{L}S^{n})\cong\left\{ \begin{array}{c}
\Lambda(a)\otimes\mathbb{Z}[u]\quad for\: n\: odd\\
(\Lambda(b)\otimes\mathbb{Z}[a,v])/(a^{2},ab,2av)\quad for\: n\: even\end{array}\right.\]
in the Morse spectral sequence we have $a\in \mathbb E^{\infty}_{0,-n}$, $b\in\mathbb E^{\infty}_{1,-2}$, 
$u\in \mathbb E^{\infty}_{1,n-2}$ and $v\in \mathbb E^{\infty}_{1,2n-3}$. 
\\
It may be interesting to notice that the energy filtration on spheres corresponds to the filtration of algebraic loops on spheres by their polynomial degree \cite{BC}, the same fact holds for complexe projective spaces.

\subsection{Morse-Serre spectral sequences}

In this section we define a fiberwise version of the preceding spectral sequence. This spectral sequence is associated to a fibration 
$$Fiber\rightarrow E\rightarrow B$$ 
where we have a Morse filtration of the base $B$ (in our case the length filtration of a loop space) and we use it to filter the total space $E$. In the classical case, where we have a smooth fiber bundle 
$$Fiber\rightarrow N\rightarrow M$$
of closed oriented compact manifolds and a Morse function $f:M\rightarrow \mathbb R$, we recover the Leray-Serre spectral sequence of the fiber bundle. 
This approach is very close in spirit to the treatment of Floer homology of families, see for example M. Hutchings'paper \cite{Hu}.

\subsubsection{Definition of the spectral sequence.}

Let $M$ be a $d$-dimensional compact oriented Riemannian manifold. Assume that the critical points of the energy function on $\mathcal LM$ are collected on compact critical manifolds. We also assume that each of these critical manifolds satisfy the Bott non-degeneracy condition.
\\
Let $p: X \to M$ be a Serre fibration (resp. a fiber bundle) over $M$ with fiber $F$ such that $F$ is a closed smooth manifold of finite dimension $f$. Let us consider the fibration (the fiber bundle) $\pi := \mathcal L p: \mathcal LX \to \mathcal LM$ with fiber $\mathcal LF$. As in the preceding sections, we filter $\mathcal LM$ by the critical values $\lambda_r$ of the energy functional. This filtration induces a filtration on $\mathcal LX$ that we denote by $\mathcal LX_{\leq \lambda_r}$ defined by the following pull-back diagram.
$$
\xymatrix{
\mathcal LX_{\leq \lambda_r} \ar[d]^{\pi} \ar@{^{(}->}[r] & \mathcal LX \ar[d]^{\pi} \\
\mathcal LM^{\leq \lambda_r} \ar@{^{(}->}[r] &\mathcal LM }
$$
We define $$\widetilde{\Sigma_r}:= \{ \gamma \in \mathcal LX / E(\pi(\gamma))=\lambda_r, \mbox{ and } dE(\pi(\gamma))=0 \} .$$ 
We remark that $\widetilde{\Sigma_r}$ is obtained from $\pi$ and $\Sigma_r$ by the following pull-back diagram
$$
\xymatrix{
\widetilde{\Sigma_r} \ar@{^{(}->}[r] \ar[d]^{\pi} &\mathcal  LX  \ar[d]^{\pi} \\
\Sigma_r \ar@{^{(}->}[r]  &\mathcal LM .} 
$$
Moreover for the sake of simplicity we suppose that these critical sets are connected.
By pulling back over $T\mathcal LM$, we define $\widetilde{\mu^{-}_p}$ in the following way. 
Let $\widetilde{T\mathcal LX}$ be the space defined by the following pull-back diagram
$$
\xymatrix{
\widetilde{T\mathcal LX} \ar[r]^{\widetilde{proj}} \ar[d]^{\widetilde{T \pi}} &\mathcal LX \ar[d]^{\pi} \\
T\mathcal LM \ar[r]^{proj} & \mathcal LM. }
$$
Where $proj$ is the canonical projection of $T\mathcal LM$.
Then, since 
$$T\mathcal LM_{\vert_{\Sigma_r}} \simeq \mu^-_r \oplus \mu^0_r \oplus \mu^+_r$$
the total space of the bundle $\widetilde{T\mathcal LM}_{\vert_{\widetilde{\Sigma_r}}} \stackrel{\widetilde{T\pi}}{\to} \widetilde{\Sigma_r}$ splits in three parts
$$\widetilde{T\mathcal LX}_{\vert_{\widetilde{\Sigma_r}}} \simeq \widetilde{\mu^-_r} \oplus \widetilde{\mu^0_r} \oplus \widetilde{\mu^+_r}.$$ 
We have the following commutative diagram
$$
\xymatrix{
\widetilde{T\mathcal LX}_{\vert_{\widetilde{\Sigma_r}}} \simeq \widetilde{\mu^-_r} \oplus \widetilde{\mu^0_r} \oplus \widetilde{\mu^+_r} \ar[r]^{\widetilde{T\pi}}  \ar[d]^{\widetilde{proj}} & T\mathcal LM_{\vert_{\Sigma_r}} \simeq \mu^-_r \oplus \mu^0_r \oplus \mu^+_r \ar[d]^{proj} \\
{{\widetilde{\Sigma_r}}} \ar[r]^{\pi} & \Sigma_r }.
$$

\begin{thm} \label{Morse Serre spec seq} Let $p:X\rightarrow M$ be a fiber bundle over a Riemannian manifold equipped with a metric such that the critical sets of the energy functional on $\mathcal LM$ satisfy the Bott non-degeneracy condition. 
Then the fiberwise energy filtration of $C_*(\mathcal LX)$ induces a spectral sequence called the Morse-Serre spectral sequence $\{E^r_{*,*}(\mathcal{MS})(\pi) \}_{r \in \mathbb{N}}$ converging to $H_*(\mathcal LX)$ denoted by
$$E^r_{p,q}(\mathcal{MS})(\pi) \Rightarrow H_{p+q}(\mathcal LX).$$
The $E^1$-page is given by $E^1_{p,q}(\mathcal{MS})(\pi) = H_{p+q}(\mathcal LX_{\leq \lambda_p},\mathcal LX_{\leq \lambda_{p-1}})$ so that $E^1_{p,q}(\mathcal{MS})(\pi)$ is isomorphic to the reduced homology $\widetilde{H}_{p+q}(Th(\widetilde{\mu^-_p}))$.
\end{thm}  

Proof. 
We have a pull-back square 
$$
\xymatrix{
\mathcal LX_{\leq \lambda_{p-1}} \ar[d]^{\pi} \ar@{^{(}->}[r] & \mathcal LX_{\leq \lambda_{p}} \ar[d]^{\pi} \\
\mathcal LM^{\leq \lambda_{p-1}} \ar@{^{(}->}[r] &\mathcal LM^{\leq \lambda_{p}} }
$$
the embedding of $\mathcal LX_{\leq \lambda_{p-1}}$ into $\mathcal LX_{\leq \lambda_{p}}$ has a normal bundle which is the pull-back of the normal bundle of the embedding of $\mathcal LM^{\leq \lambda_{p-1}}$ into $\mathcal LM^{\leq \lambda_{p}}$. We recall that the space $\mathcal LM^{\leq \lambda_{p}}/\mathcal LM^{\leq \lambda_{p-1}}$ is homotopy equivalent to the Thom space   $Th(\mu^-_p)$ we deduce that   $\mathcal LX_{\leq \lambda_p} /\mathcal LX_{\leq \lambda_{p-1}}$ is homotopically equivalent to the Thom space $Th(\widetilde{\mu^-_p})$.
\\
\rightline{$\square$}

\subsubsection{Compatibility with the loop product.}
We assume that $M$ satisfies the condition (Cl).
\\
The embedding $\delta_X: X \hookrightarrow X \times X, x \mapsto (x,x)$ factorizes by the map $\delta^1_X$ defined by the following pull back diagram
$$
\xymatrix{
X \times_M X \ar@{^{(}->}[r]^{\delta^1_X} \ar[d]^{p \times_M p} & X \times X \ar[d]^{p \times p} \\
M \ar@{^{(}->}[r]^{\delta_M} & M \times M }
$$
and by the canonical map $\delta^2_X: X \hookrightarrow X \times_M X$ induced by the diagonal map on $X$ and the projection on $M$. Thus, we have $\delta_X = \delta^1_X \circ \delta^2_X$. We denote by $\widetilde{\delta^i_X}$ the pull-back of $\delta^i_X$ over the free loop fibration, $i=1,2$. We remark that $\delta^1_X$ is a smooth embedding of codimension $dim(F)$, we denote it by $f$.

Let us recall that $F_pC_{p+q}(\mathcal LX) = C_{p+q}(\mathcal LX_{\leq \lambda_p})$ then the loop product is given by the following composition map.
$$
\xymatrix{
C_{p+q}(\mathcal LX_{\leq \lambda_p}) \otimes C_{p'+q'}(\mathcal LX_{\leq \lambda_{p'}}) \ar[d]^{\times} \\
C_{p+q+p'+q'}(\mathcal LX_{\leq \lambda_p} \times \mathcal LX_{\leq \lambda_{p'}}) \ar[d]^{\widetilde{\delta^1_X}!} \\
C_{p+q+p'+q'-d}(\mathcal LX_{\leq \lambda_p} \times_M \mathcal LX_{\leq \lambda_{p'}}) \ar[d]^{\widetilde{\delta^2_X}!} \\
C_{p+q+p'+q'-d-f}(\mathcal LX_{\leq \lambda_p} \times_X \mathcal LX_{\leq \lambda_{p'}}) \ar[d]^{comp_{X_*}} \\
C_{p+q+p'+q'-d-f}(\mathcal LX_{\leq \lambda_p+\lambda_{p'}}) = C_{p+q+p'+q'-d-f}(\mathcal LX_{\leq \lambda_{p+p'}}) }.
$$

This gives a multiplicative structure to the $(0, d+f)$-regraded Morse Serre spectral sequence. We have the following proposition 
\begin{proposition} \label{mult Morse Serre}
Let $M$ be a Riemannian compact $d$-dimensional manifold that satisfies condition (Cl). Let $\pi: X \to M$ be a fibration over $M$.
Then, the loop product induces a multiplicative structure on the regraded Morse Serre spectral sequence $\mathbb{E}^r_{p,q}(\mathcal{MS})(\pi) := E^r_{p, q+ d + f}(\mathcal{MS})(\pi)$.
\end{proposition}

\subsubsection{The first page} Now let us focus on the $E^1$-page of the Morse Serre spectral sequence.
In order to describe this first page and its multiplicative structure we will introduce a bigraded algebra 
$\mathcal A^F_{*,*}$, this bigraded algebra is the fibrewise analogue of the bigraded algebra $\mathcal A_{*,*}$ introduced in the preceding sections.
\\
Let $\widetilde{M}$ and $\widetilde{UM}$ the fiberwise analogues of the spaces $M$ and $UM$ : 
\\
a) the space $\widetilde{M}$ is defined as the pull-back of the inclusion of constant loops $M\subset \mathcal LM$ along the fibration $\pi:\mathcal LX\rightarrow LM$. We thus have a fibration $c:\widetilde{M}\rightarrow M$.
\\
b) the space $\widetilde{UM}$ is defined as the pull-back of the map $geod:UM\rightarrow \mathcal LM$ which associates to a unit tangent vector $u$ the unique closed geodesic such that $\dot{\gamma}(0)=u$. We also have a fibration $b: \widetilde{UM} \to UM$.
\\
We define a multiplicative structure on $\mathbb{H}_*(\widetilde{UM}) := H_{*+d + d-1 +f}(\widetilde{UM}) $ that mixes the intersection product on $UM$ and the loop product on the fiber $LF$ of $\pi$.
This product  is associated to the following commutative diagram
$$
\xymatrix{
\mathcal LF \times \mathcal LF \ar[r] & \widetilde{UM} \times \widetilde{UM} \ar[r]^{b \times b} & UM \times UM \\
\mathcal LF \times \mathcal LF \ar[r] \ar[u]^{=} & \widetilde{UM} \times_{UM} \widetilde{UM} \ar[r]^{b \times_{UM} b} \ar@{^{(}->}[u]^{\widetilde{\delta_{UM}}} & UM \ar@{^{(}->}[u]^{\delta_{UM}} \\
\mathcal LF \times_F \mathcal LF \ar[d]^{comp_F} \ar[r] \ar@{^{(}->}[u]^{\widetilde{\delta_F}} & \widetilde{UM}^{\infty} \ar[d]^{comp_{\widetilde{UM}^{\infty}}} \ar[r]^{b} \ar@{^{(}->}[u]^{\hat{\delta_F}} & UM \ar[u]^{=} \ar[d]^{=}\\
\mathcal LF \ar[r] & \widetilde{UM} \ar[r]^{b}  & UM }
$$
where $\widetilde{UM}^{\infty}$ denotes the composable loops of $\widetilde{UM}$. The product
 is given by the composition
$$comp_{{\widetilde{UM}^{\infty}}_*} \circ \hat{\delta_F}! \circ \widetilde{\delta_{UM}}! \circ \times.$$ 
\\
In the same way we a define a multiplicative structure for the homology of $\widetilde{M}$ by considering it as the total space of the fibration $c$. 
This provides $\mathbb{H}_*(\widetilde{M}):= H_{*+d+f}(\widetilde{M})$ with a structure of algebra.
\\
Furthermore, there is a structure of $\mathbb{H}_*(\widetilde{M})$-module on $\mathbb{H}_*(\widetilde{UM})$. This structure comes from the following commutative diagram
$$
\xymatrix{
\mathcal LF \times \mathcal LF \ar[r] & \widetilde{M} \times \widetilde{UM} \ar[r]^{c \times b} & M \times UM \\
\mathcal LF \times \mathcal LF \ar[r] \ar[u]^{=} & \widetilde{M} \times_M \widetilde{UM} \ar[r] \ar@{^{(}->}[u]^{\widetilde{\hat{\delta_M}}} & UM \ar@{^{(}->}[u]^{\hat{\delta_M}} \\
\mathcal LF \times_F \mathcal LF \ar[d]^{comp_F} \ar[r] \ar@{^{(}->}[u]^{\widetilde{\delta_F}} & \widetilde{UM}^{\infty} \ar[r] \ar@{^{(}->}[u]^{\underline{\widetilde{\delta_F}}} \ar[d]^{comp_{\widetilde{UM}^{\infty}}} & UM \ar[u]^{=} \ar[d]^{=}\\
\mathcal LF \ar[r] & \widetilde{UM} \ar[r]^{b}  & UM }.
$$
where $\hat{\delta_M}$ denotes the following inclusion
$$\hat{\delta_M}: UM \hookrightarrow M \times UM $$
$$ x \mapsto (proj(x), x)$$
with $proj: UM \to M$ is the canonical projection. The module structure $$\widetilde{m}: \mathbb{H}_*(\widetilde{M}) \otimes \mathbb{H}_*(\widetilde{UM}) \to \mathbb{H}_*(\widetilde{UM})$$ is given by
$$\widetilde{m} = comp_{{\widetilde{UM}^{\infty}}_*} \circ \underline{\widetilde{\delta_F}}! \circ {\widetilde{\hat{\delta_M}}}! \circ \times.$$

\begin{definition}
 These products on $\mathbb{H}_*(\widetilde{UM})$ and on $\mathbb{H}_*(\widetilde{M})$ are called the fiberwise intersection loop product. The module structure $\widetilde{m}$ of $\mathbb{H}_*(\widetilde{M})$ on $\mathbb{H}_*(\widetilde{UM})$ is called the fiberwise loop product module structure. The bigraded algebra $\mathcal A^F_{*,*}$ is defined by
\\
1) $\mathcal A^F_{0,*}=\mathbb{H}_*(\widetilde{M})$
\\
2) $\mathcal A^F_{p,q}=\mathbb{H}_{q-p \alpha_1}(\widetilde{UM})<T^p>$ with $T$ of bidegree $(1,\alpha_1+d-2)$ for $p>0$
\\
the multiplicative structure is given by the fiberwise loop product.

\end{definition}

Let us give the main Theorem of this section. This theorem will be one of the main computational tool of this paper.

\begin{thm} \label{morseserreE1}
Let $M$ be a $d$-dimensional compact oriented Riemannian manifold that satisfies condition (Cl). Let $p: X \to M$ be a fibration over $M$.
The $E^1$-page of the regraded Morse Serre spectral sequence 
$\mathbb{E}^r_{p,q}(\mathcal{MS})(\pi) := E^r_{p, q+ d + f}(\mathcal{MS})(\pi)$ is isomorphic as a bigraded algebra to 
$\mathcal A^F_{*,*}$.
\end{thm}

Proof.
For $r \geq 1$, we define $$\widetilde{\phi_r} :H_*(\mathcal LX_{\leq \lambda_r},\mathcal  LX_{\leq \lambda_{r-1}}) \to H_{*- \alpha_r}(\widetilde{UM})$$ by composition of the Thom isomorphism:
$$\mathcal{T}: H_*(\mathcal LX_{\leq \lambda_r},\mathcal  LX_{\leq \lambda_{r-1}}) \to H_{*- \alpha_r}(\widetilde{\Sigma_r})$$ 
with the identification of $H_{*- \alpha_r}(\widetilde{\Sigma_r})$ with $H_{*- \alpha_r}(\widetilde{UM})$. Let us denote
$\mathcal LX_{\leq \lambda_r}$ by $\Lambda_r$, the proof follows from the following commutative diagram
$$
\xymatrix{
H_p(\Lambda_r,\Lambda_{r-1}) \otimes H_{p'}(\Lambda_{r'},\Lambda_{r'-1}) \ar[r]^{\widetilde{\phi_r} \otimes \widetilde{\phi_r}} \ar[d]^{\times} & H_{p- \alpha_r}(\widetilde{UM}) \otimes H_{p'- \alpha_{r'}}(\widetilde{UM}) \ar[d]^{\times} \\
H_{p+p'}(\Lambda_r \times \Lambda_{r'},\Lambda_r\times\Lambda_{r-1}\cup \Lambda_{r-1}\times \Lambda_{r'}) \ar[r]^{\widetilde{\phi_r} \times \widetilde{\phi_r}} \ar[d]^{\widetilde{\delta_X}_!} & H_{p + p'- \alpha_r - \alpha_{r'}}(\widetilde{UM} \times \widetilde{UM} ) \ar[d]^{\hat{\delta_F}_! \circ \widetilde{\delta_{UM}}_!} \\
H_{p+p'-d-f}(\Lambda_r \times_X \Lambda_{r'},\Lambda_r\times_X\Lambda_{r-1}\cup \Lambda_{r-1}\times_X \Lambda_{r'}) \ar[r]^{\widetilde{\phi_r} \times_X \widetilde{\phi_r}} \ar[d]^{{comp_X}_*} & H_{p + p'- \alpha_r - \alpha_{r'} -d - (d-1) -f}(\widetilde{UM}^{\infty} ) \ar[d]^{{comp_{\widetilde{UM}^{\infty}}}_*} \\
H_{p+p'-d-f}(\Lambda_{r+r'}) \ar[r]^{\widetilde{\phi_{r +r'}}} & H_{p + p'- \alpha_r - \alpha_{r'} -d - (d-1) -f}(\widetilde{UM}) }.
$$
The last line can be rewritten 
$$
H_{p+p'-d-f}(\Lambda_{r+r'},\Lambda_{r+r'-1}) \to H_{p + p'- \alpha_{r+r'} -d - (d-1) -f}(\widetilde{UM}).$$
\\
For $p=0$, we notice that $\mathbb{H}_*(\mathcal LX_{\leq 0},\mathcal  LX_{\leq -1})$ with the loop product is in fact $\mathbb{H}_*(\widetilde{M})$ with the fiberwise intersection loop product. This allows the identification of the $0^{th}$ column 
$$\mathbb{E}^1_{0,q}(\mathcal{MS})(\pi)= \mathbb{H}_q(\widetilde{M}).$$
\rightline{$\square$}

\section{The $0$-th column of the spectral sequence}

In this section we begin the computation of $\mathbb H_*(Imm(S^1,S^n))$. By using Hirsch-Smale theory we replace the space $Imm(S^1,S^n)$ by the loop space $\mathcal LUS^n$ and we consider the fiber bundle
$$\mathcal LS^{n-1}\rightarrow \mathcal LUS^n\rightarrow \mathcal LS^n$$
associated to the unit tangent bundle
$$S^{n-1}\rightarrow US^n\rightarrow S^n.$$
We filter the space $\mathcal LUS^n$ by the length filtration of the base $\mathcal LS^n$ and compute the associated spectral sequence. We notice that the $0$-th column of the $E^1$-term of this spectral sequence is the homology of the space of vertical loops of $\mathcal LUS^n$. Vertical loops are the loops of $US^n$ constant on $S^n$. We denote this space by $\mathcal LUS^n_{v}$

\subsection{The even case}

\begin{proposition} \label{E^0 pair}
The algebra $\mathbb H_*(\mathcal LUS^{2n}_v)$ is isomorphic to $$\mathbb{Z}[a, b, c ]/(b^2, 2b, a^2, ab)$$
with $deg(a)=-4n+1$, $deg(b)=-2n$ and $deg(c)=2n-2$.
\end{proposition}

\noindent{Proof.}

By Theorem \ref{morseserreE1} we have the isomorphism $\mathbb{E}^1_{0,*}(\mathcal{MS})(\pi)\cong \mathbb{H}_*(\widetilde{S^{2n}})$.
Let us determine the $(2n, 2n+1)$-regraded Serre spectral sequence $\mathbb{E}_{*,*}^*(c)$ associated to the fibration 
$c: \widetilde{S^{2n}} \to S^{2n}$ with fiber $\mathcal LS^{2n-1}$.
This spectral sequence is multiplicative and 
$$\mathbb{E}^2_{*,*}(c) = \mathbb{H}_*(S^{2n}; \mathbb{H}_*(\mathcal LS^{2n-1})) = \mathbb{H}_*(S^{2n}) \otimes  \mathbb{H}_*(\mathcal LS^{2n-1})$$
since $\mathbb{H}_*(S^{2n})$ has no torsion.
Then, $\mathbb{E}^2_{*,*}(c)$ is isomorphic to $$\mathbb{Z}[x_{-2n}]/x_{-2n}^2 \otimes \mathbb{Z}[u_{2n-2}, y_{-2n+1}]/y_{-2n+1}^2.$$
On the following diagram and on the others of this paper, the symbols X represent the abelian group $\mathbb{Z}$ and the symbols O the field $\mathbb{Z}/2\mathbb{Z}$. The generators of algebra are represented by squares.

\begin{center}
\includegraphics[width=11cm]{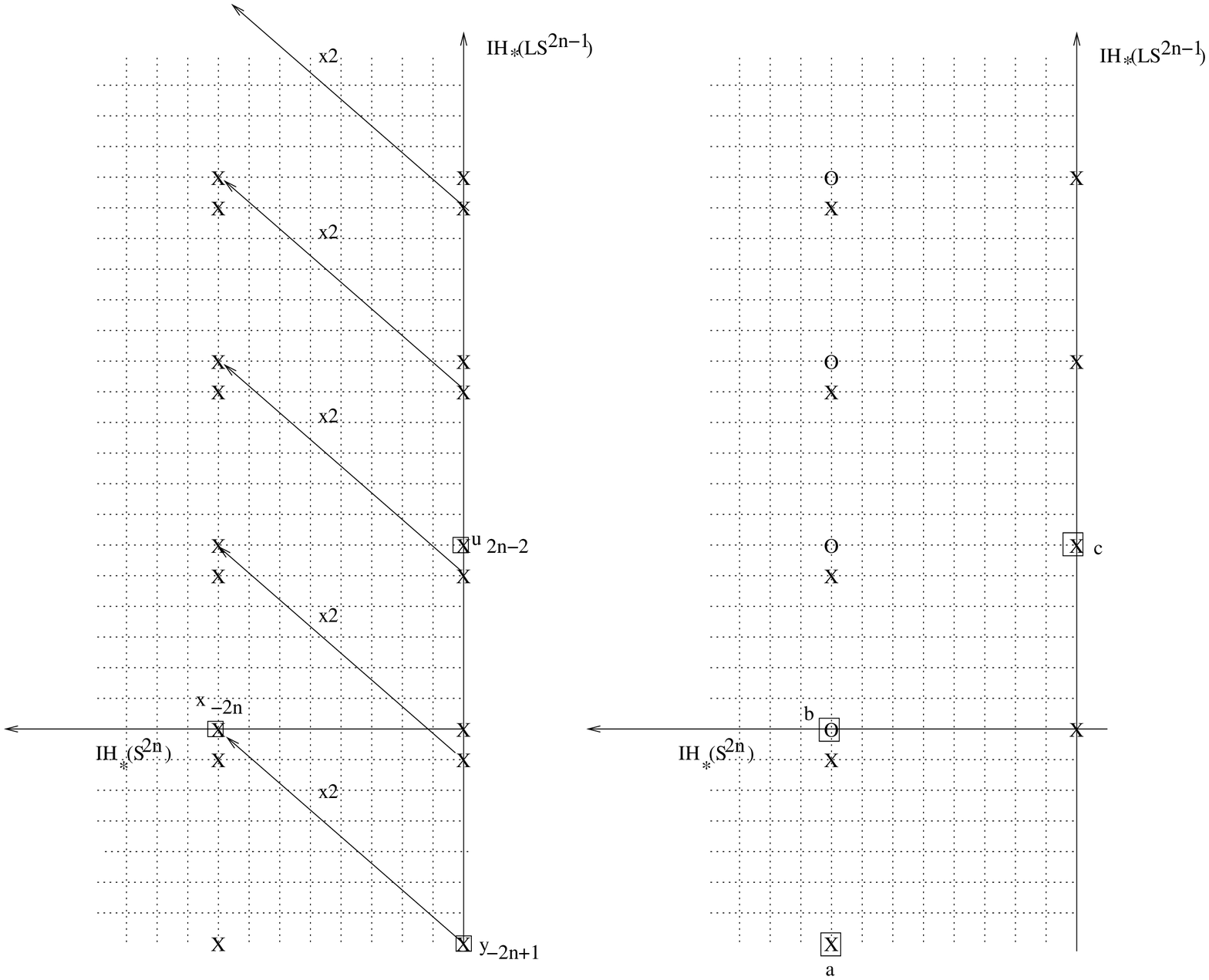} 
\end{center}

Since this spectral sequence is a spectral sequence of algebra, we only need to compute the differentials on the generators. For reasons of degree, the only possible non zero differential is $d_{2n}(y_{-2n+1})$. The map $ev(0):\mathcal LUS^{2n}_{\leq 0} \to US^{2n}$ has a section so that $\mathbb{H}_{-2n}(US^{2n})$ is isomorphic to $\mathbb{H}_{-2n}(\widetilde{S^{2n}})$ namely $\mathbb{Z}/2\mathbb{Z}$. This implies that $d_{2n}(y_{-2n+1}) = 2x_{-2n}$.

For degree reasons and by position in the filtration, there is neither linear nor multiplicative extension issues so that $$\mathbb{H}_*(\widetilde{S^{2n}}) \cong \mathbb{Z}[w_{-4n+1}, x_{-2n}, u_{2n-2}]/(w_{-4n+1}^2, w_{-4n+1}x_{-2n}, x_{-2n}^2,2x_{-2n})$$
where $w_{-4n+1} := x_{-2n}y_{-2n+1}$.
We denote by $a$,$b$ and $c$ respectively the image of $w_{-4n+1}$, $x_{-2n}$ and $u_{2n-2}$  under the isomorphism of graded algebra
from $\mathbb{H}_*(\widetilde{S^{2n}})$ to $\mathbb{E}^1_{0,*}(\mathcal{MS})(\pi)$. The degrees are the same namely $deg(a)=-4n+1$, $deg(b)=-2n$ and $deg(c)=2n-2$.

\rightline{$\square$}

\subsection{The odd case}

\begin{proposition} \label{E^0 impair}
The algebra $\mathbb H_*(\mathcal LUS^{2n+1}_v)$ is isomorphic to $$\mathbb{Z}[\alpha, \beta, \gamma, \delta ]/(\alpha^2, \gamma^2, \beta^2, 2 \beta \delta)$$
with $deg(\alpha)=-2n-1$, $deg(\beta)=-2n$, $deg(\gamma)=-1$ and $deg(\delta)=4n-2$.
\end{proposition}

Proof.
As in the preceding subsection, we compute $\mathbb{E}^{1}_{0,*}(\mathcal{MS})(\pi)$  by computing the $(2n+1, 2n)$-regraded  Serre spectral sequence associated to the fibration $c: \widetilde{S^{2n+1}} \to S^{2n+1}$.

\begin{center}
\includegraphics[width=4cm]{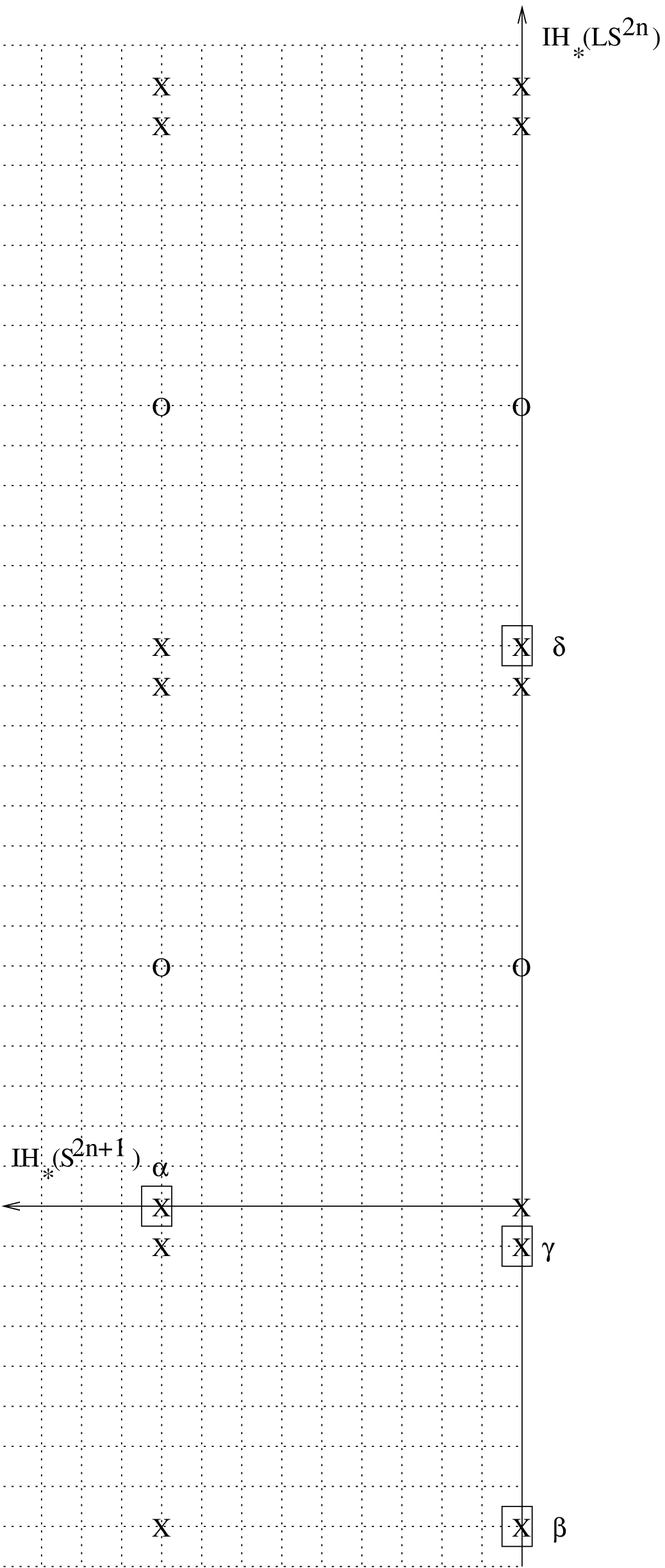} 
\end{center}

In this case, the spectral sequence collapses at the $E^2$-level because the fibration $ev(o): \mathcal L US^{2n}_{\leq 
0} \to US^{2n}$ admits a section and so $\mathbb{H}_{-2n-1}(\widetilde{S^{2n+1}}) \simeq \mathbb{Z}$. For reasons of degree, there are no extension issues.

\rightline{$\square$}



\section{Computation of $\mathbb{H}_*(\mathcal LUS^{2n})$, $n \geq 2$.} \label{LUS2n}

\noindent{\bf Strategy of the proof.} 

(1) We begin by computing $\mathbb{H}_*(\widetilde{US^{2n}})$ in subsection \ref{mathbb}, we will get $\mathbb{E}^1_{1,*}(\mathcal{MS})(\pi)$. We do it by computing the Serre spectral sequence associated to the fibration $b:\widetilde{US^{2n}} \to US^{2n}$. 

(2) Then, we compute in subsection \ref{paire module structure} the $\mathbb{H}_*(\widetilde{S^{2n}})$-module structure on $\mathbb{H}_*(\widetilde{US^{2n}})$. This will give the multiplicative structure of the first page of the Morse-Serre spectral sequence $\mathbb{E}^1_{*,*}(\mathcal{MS})(\pi)$. 

(3) Next, in \ref{paire CJY}, we compute the $E^{\infty}$-page of the Cohen-Jones-Yan spectral sequence (that is also the $E^2$-page).
This gives all the differentials on $\mathbb{E}^1_{*,*}(\mathcal{MS})(\pi)$ and leads to $\mathbb{E}^{\infty}_{*,*}(\mathcal{MS})(\pi)$.

(4) Finally, we compare the $E^{\infty}$-page of the two spectral sequences and we solve the linear extension issues and the multiplicative extension issues.

The different stages of the proof are illustrated by pictures representing the different spectral sequences for $n=4$ namely in the case $US^8$.

\subsection{The columns of the Morse-Serre spectral sequence: $\mathbb{E}^1_{p,*}(\mathcal{MS})(\pi)$, $p \geq 1$.} \label{mathbb}
In this subsection, we compute the $p^{th}$ columns of the Morse Serre spectral sequence of $US^{2n}$ for $p \geq 1$.

\begin{proposition} \label{E^1}
The subalgebra $\mathbb{E}^1_{p,*}(\mathcal{MS})(\pi)$, $p \geq 1$, is isomorphic to the algebra
$$\mathbb{Z}[i, j, k, l ]/(i^2, j^2, 2j, ij, k^2)[T]_{\geq 1}$$
with $deg(i)=-4n+1$, $deg(j)=-2n$, $deg(k)=-2n+1$, $deg(l)=2n-2$ and $bideg(T)= (1,4n-3)$.
\end{proposition}

\noindent{Proof.}

It suffices to compute $\mathbb{E}^1_{1,*}(\mathcal{MS})(\pi)$, the other columns are isomorphic up to a shifting and are given by multiplication by $T$.

We begin by applying Theorem \ref{morseserreE1}. We first need to compute $\mathbb{H}_{*}(\widetilde{US^{2n-1}})$. Let us prove that the algebra $\mathbb{H}_{*}(\widetilde{US^{2n-1}})$ is isomorphic to $$\mathbb{Z}[x_{-4n+1}, z_{-2n}, u_{2n-2}, y_{-2n+1}]/(x_{-4n+1}^2, 2z_{-2n}, z_{-2n}^2, y_{-2n+1}^2).$$
\\
In order to proceed let us compute the $(4n-1, 2n-1)$-regraded Serre spectral sequence associated to the fibration $b:\widetilde{US^{2n}} \to US^{2n}$ with fiber $\mathcal LS^{2n-1}$. We have $$\mathbb{E}^2_{*,*}(b) \cong \mathbb{H}_*(US^{2n}; \mathbb{H}_*(LS^{2n-1})).$$ 
Since $\mathbb{H}_*(\mathcal LS^{2n-1})$ has no torsion, we also have $$\mathbb{H}_*(US^{2n}; \mathbb{H}_*(\mathcal LS^{2n-1})) \cong \mathbb{H}_*(US^{2n}) \otimes  \mathbb{H}_*(\mathcal LS^{2n-1})$$ namely this algebra is isomorphic to 
$$\mathbb{Z}[x_{-4n+1}, z_{-2n}, u_{2n-2}, y_{-2n+1}]/(x_{-4n+1}^2, 2z_{-2n}, z_{-2n}^2, y_{-2n+1}^2).$$ 

On the following spectral sequence and on the others of the article, the symbol O represents the abelian group $\mathbb{Z}/2\mathbb{Z}$.

\begin{center}
\includegraphics[width=6cm]{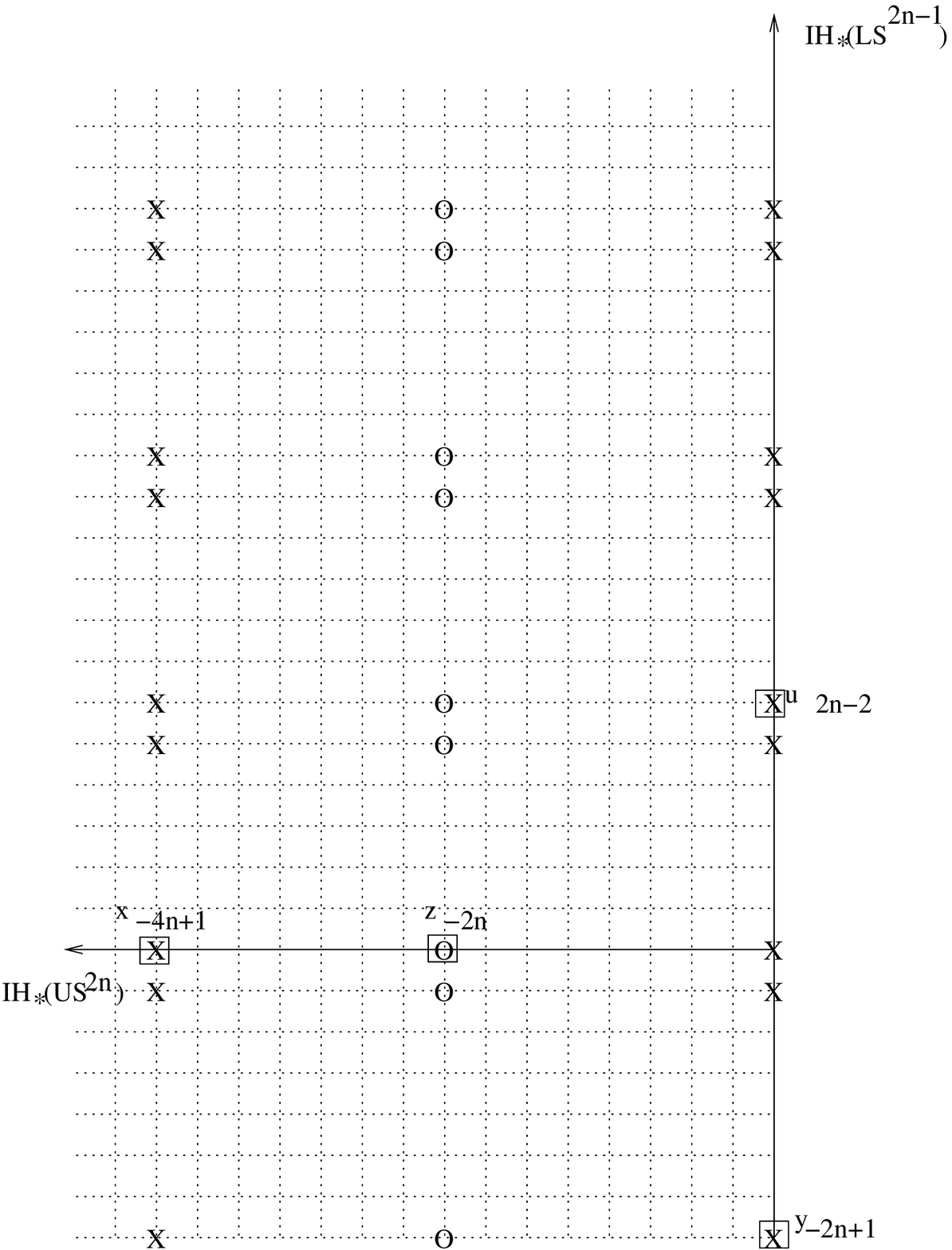} 
\end{center}

Let us prove that this spectral sequence collapses at the $E^2$-level. We first notice that the fibration $\mathcal \pi :\widetilde{\Sigma_1} \to \Sigma_1$ admits a section. This section is the map
$$s: \Sigma_1 \to \widetilde{\Sigma_1}$$
$$ \gamma \mapsto (\gamma, \frac{\dot{\gamma} }{\Vert \dot{\gamma} \Vert}).$$
By the bundle isomorphism $\psi$ of proposition \ref{mult Morse Serre}, this proves that $b:\widetilde{US^{2n}} \to US^{2n}$ admits a section.
Then, for degree reasons, the other differentials vanish.  To achieve the proof we need to solve the extension issues. The only extension issue is the product $z_{-2n} y_{-2n+1} $ which is defined modulo $\mathbb{Z}<x_{-4n+1}>$.  But $z_{-2n}$ is of $2$-torsion, then $z_{-2n} y_{-2n+1}$ is also of $2$-torsion. Then the component of $z_{-2n} y_{-2n+1} $ is zero. The same argument works for all products $z_{-2n} y_{-2n+1} u_{2n-2}^k$, $k \in \mathbb{N}$.
\rightline{$\square$}

We represent the first page of the Morse-Serre spectral sequence $\mathbb{E}^1_{*,*}(\mathcal{MS})(\pi)$ on the following diagram.

\begin{center}
\includegraphics[width=10cm]{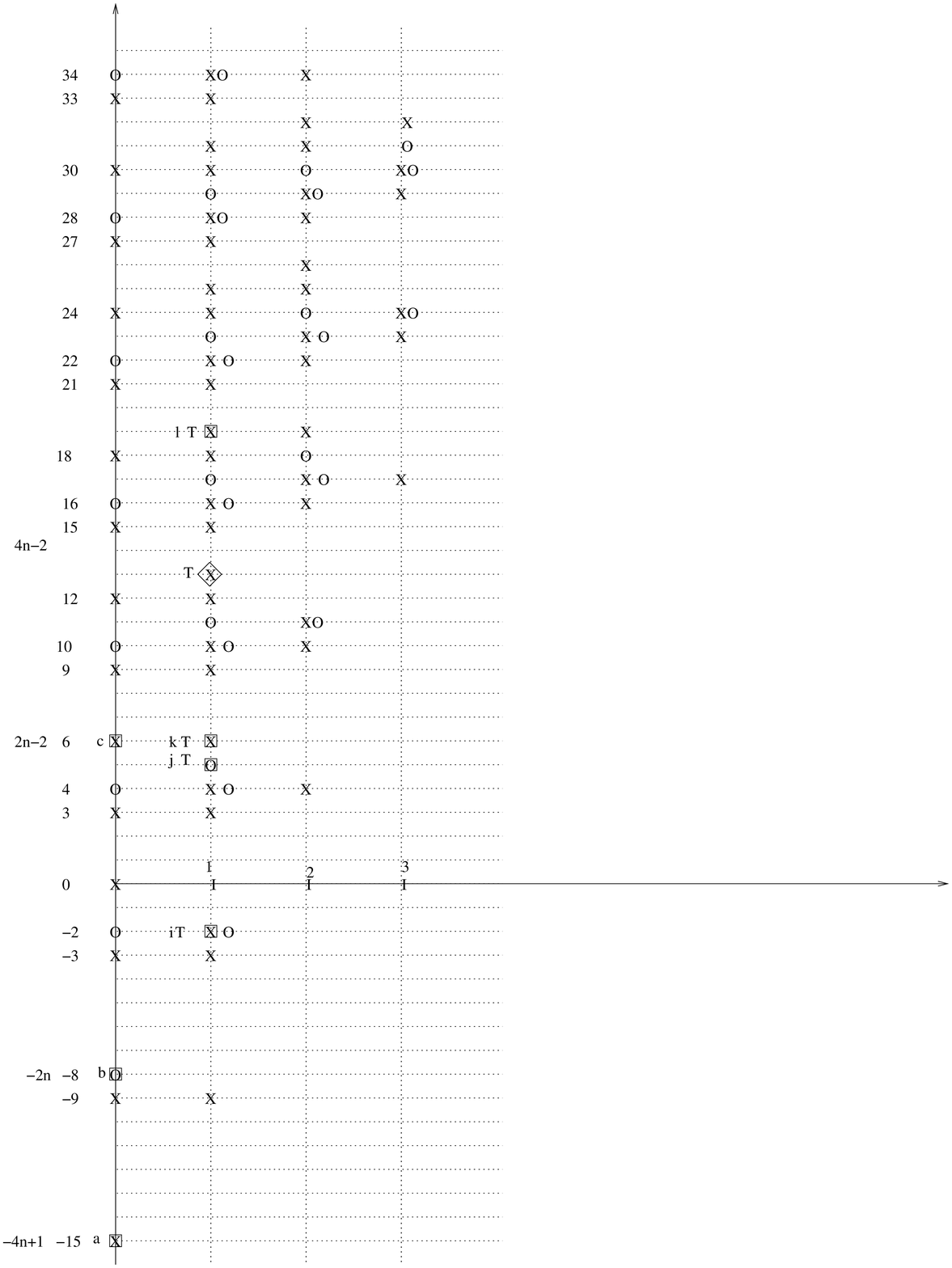} 
\end{center}

\subsection{The $\mathbb{E}^1_{0,*}(\mathcal{MS})(\pi)$-module structure on $\mathbb{E}^1_{1,*}(\mathcal{MS})(\pi)$.} \label{paire module structure}

\begin{proposition} \label{module}
The structure of $\mathbb{H}_*(\widetilde{S^{2n}})$-module on $\mathbb{H}_*(\widetilde{US^{2n}})$ is entirely described by the following products: $c.i= li$, $c.j=lj$, $c.k=lk$, $c.l=l^2$, $b.k=jk$, $b.l=lj$, $a.l=klj$. The other products are zero.
\end{proposition}

\noindent{Proof.}

It follows from the commutative diagram of section 3.2 that we have a $\mathbb{E}^*_{*,*}(c)$-module structure on $\mathbb{E}^*_{*,*}(b)$. Then, we will compute this module structure using these spectral sequences. At the $E^2$-level, this module structure is given by the Chas and Sullivan loop product on the fiber and by the $\mathbb{H}_*(S^{2n})$-module structure on $\mathbb{H}_*(US^{2n})$ induced by $\hat{\delta_M}$ (see section 3.2).
First, let us recall that 
$$\mathbb{H}_*(S^{2n}) \cong \mathbb{Z}[k_{-2n}]/(k_{-2n}^2)$$
and that
$$\mathbb{H}_*(US^{2n}) \cong \mathbb{Z}[g_{-4n+1}, h_{-2n}]/(g_{-4n+1}^2, g_{-4n+1}h_{-2n}, 2h_{-2n}, h_{-2n}^2).$$

\begin{lem} \label{l}
 The $\mathbb{H}_*(S^{2n})$-module structure on $\mathbb{H}_*(US^{2n})$ is given by the following products:
$k_{-2n}.g_{-4n+1}=0$, $k_{-2n}.h_{-2n}=0$ and $k_{-2n}.1_{US^{2n}}=h_{-2n}$.
\end{lem}

Proof.

The two first products vanish for degree reasons. For the last product, let us consider the canonical fibration $f:US^{2n} \to S^{2n}$.
By Theorem 2 of \cite{L}, the following pullback diagram
$$
\xymatrix{
US^{2n} \ar@{^{(}->}[r] \ar[d] & S^{2n} \times US^{2n} \ar[d]^{id \times f} \\
S^{2n} \ar@{^{(}->}[r]^{\delta_{S^{2n}}} & S^{2n} \times S^{2n} }
$$
induces a structure of $\mathbb{H}_*(S^{2n})$-module on the shifted Serre spectral sequence of the fibration $f:US^{2n} \to S^{2n}$.
This structure is easy to compute. There is no extension issue for degree reasons.
This proves Lemma \ref{l}.
\rightline{$\square$}

Now, let us consider the two shifted multiplicative Serre spectral sequences associated to the fibrations $b:\widetilde{US^{2n}} \to US^{2n}$ and $c:\widetilde{S^{2n}} \to S^{2n}$.
At the $E^2$-level, the module structure is given by the module structure of the preceding Lemma \ref{l} on the base and by the Chas and Sullivan loop product on the fiber. Then, we have:

\begin{enumerate}
\item $b.u_{2n-2}= z_{-2n}.u_{2n-2}$ and $b.y_{-2n+1}= z_{-2n}y_{-2n+1}$
\item $a. u_{2n-2}= z_{-2n}u_{2n-2}y_{-2n+1}$
\item $c.u_{2n-2}= u_{2n-2}^2$, $c.y_{-2n+1}= u_{2n-2}y_{-2n+1}$, $c.z_{-2n}=u_{2n-2}z_{-2n}$ and $c.x_{-4n+1}= u_{2n-2}x_{-4n+1}$.
\end{enumerate}

The other products vanish. For degree reasons, the only possible extension issue may be $c.y_{-2n+1}z_{-2n}$. But since this product lifts to an element of $2$-torsion, there is no extension issue.

\rightline{$\square$}

\subsection{The multiplicative structure of $\mathbb{E}^1_{*,*}(\mathcal{MS})(\pi)$.}
This subsection is the application of Theorem \ref{morseserreE1}. Since we know all the columns of the Morse-Serre spectral sequence and its multiplicative structure, we have the following result.

\begin{proposition} \label{E1 paire}
The algebra $\mathbb{E}^1_{*,*}(\mathcal{MS})(\pi)$ is isomorphic to
$$\mathbb{Z}[a, b, c ]/(b^2, 2b, a^2, ab) \oplus \mathbb{Z}[i, j, k, l]/(i^2, 2j, j^2, k^2)[T].$$
with $deg(a)=-4n+1$, $deg(b)=-2n$,$deg(c)=2n-2$ and $deg(i)=-4n+1$, $deg(j)=-2n$, $deg(k)=-2n+1$,$deg(l)=2n-2$, $bideg(T)=(1,4n-3)$.
The multiplication between an element of the $0^{th}$ column and an element of another column is given by $c.i= l.i$, $c.j=l.j$, $c.k=l.k$, $c.l=l^2$, $b.k=j.k$, $b.l=l.j$, $a.l=k.l.j$.
\end{proposition}

This can be illustrated by the following diagram. On this diagram, we have represented the generators $a$, $b$, $c$ of the $0$-th column and the generators $i$, $j$, $k$, $l$ of the first column by squares. We also put a diamond $\diamond$ around the generator $T$. We recall that this generator is of bidegree $(1, 4n-3)$ so that its total degree is $4n-2=\alpha_1$.

\begin{center}
\includegraphics[width=13cm]{morseserreUS8E1.eps} 
\end{center}

\subsection{Comparison with the Cohen-Jones-Yan spectral sequence.} \label{paire CJY}

Let us determine the differentials of the Morse-Serre spectral sequence. Let us recall the $E^{\infty}$-page of the Cohen-Jones-Yan spectral sequence computed in \cite{L}.

\begin{center}
\includegraphics[width=7cm]{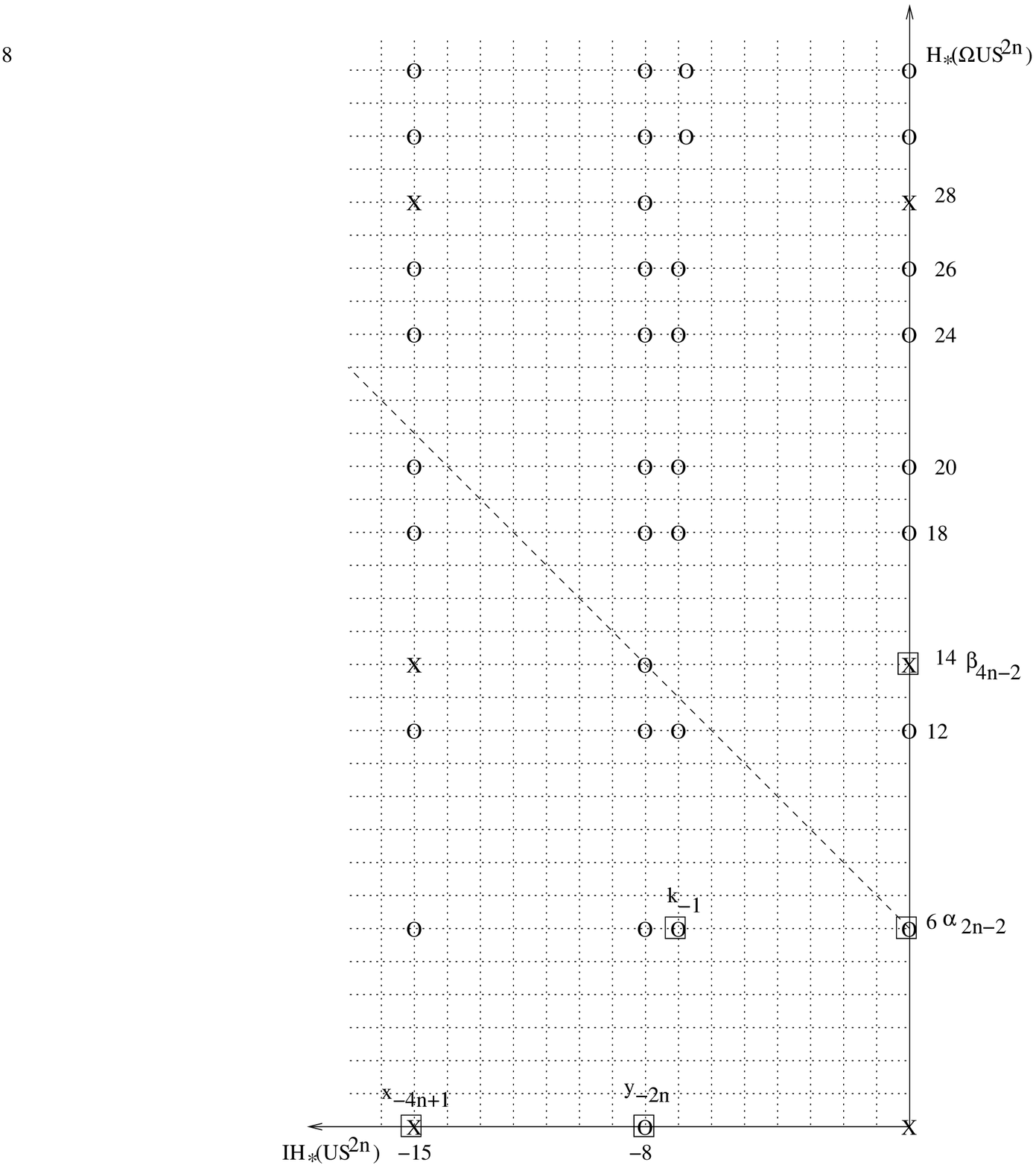} 
\end{center}

The following proposition describes the multiplicative structure of this spectral sequence.

\begin{proposition} \label{multCJY}
The algebra $\mathbb{E}^{\infty}_{*,*}[ev(0)]$ is isomorphic to

$$\mathbb{Z}[x_{-4n+1}, y_{-2n}, \alpha_{2n-2}, \beta_{4n-2}, k_{-1}]/\mathcal I$$
with $\mathcal I$ be the ideal generated by 
$$\left\lbrace x_{-4n+1}^2,x_{-4n+1}y_{-2n}, x_{-4n+1}k_{-1}, 2y_{-2n}, y_{-2n}^2, y_{-2n}k_{-1}-x_{-4n+1}\alpha_{2n-2}, 2k_{-1}, 2\alpha_{2n-2} \right\rbrace .$$

\end{proposition} 

\noindent{Proof.}
As the Cohen-Jones-Yan spectral sequence collapses at the $E^2$-term, we just have to give a description by generators and relations of the graded algebra $\mathbb{H}_*(US^{2n} ; H_*(\Omega US^{2}))$. 
\\
We first notice that this algebra is generated by the sub-algebra $\mathbb{H}_*(US^{2n}) \otimes H_*(\Omega US^{2n})$ and an element $k_{-1}\in \mathbb H_{-2n+1}(US^{2n},H_{2n-2}(\Omega US^{2n}))=\mathbb{Z} / 2\mathbb{Z}$. We recall that 
$$\mathbb{H}_*(US^{2n}) \otimes H_*(\Omega US^{2n})\cong \mathbb{Z}[x_{-4n+1}, y_{-2n}, \alpha_{2n-2}, \beta_{4n-2}]/\mathcal J$$
where the ideal $\mathcal J$ is generated by 
$$\left\lbrace x_{-4n+1}^2,x_{-4n+1}y_{-2n}, 2y_{-2n}, y_{-2n}^2, 2\alpha_{2n-2} \right\rbrace .$$
Thus we are left with the computations of the relations involving $k_{-1}$ and the sub-algebra $\mathbb{H}_*(US^{2n}) \otimes H_*(\Omega US^{2n})$. As $k_{-1}$ is of $2$-torsion we compute the Cohen-Jones-Yan spectral sequence for the coefficient ring $\mathbb{Z} / 2\mathbb{Z}$ that we denote by $\mathbb{E}^{r}_{*,*}[ev(0);\mathbb{Z} / 2\mathbb{Z}]$.
This spectral sequence also collapses at the $E^2$-term and as we work over a field we have $$\mathbb{E}^{\infty}_{*,*}[ev(0);\mathbb{Z} / 2\mathbb{Z}] = \mathbb{H}_*(US^{2n};\mathbb{Z} / 2\mathbb{Z} ) \otimes H_*(\Omega US^{2n}; \mathbb{Z} / 2\mathbb{Z}).$$
And the algebra morphism
$$\mathbb{E}^{\infty}_{*,*}[ev(0);\mathbb{Z}]\rightarrow\mathbb{E}^{\infty}_{*,*}[ev(0);\mathbb{Z} / 2\mathbb{Z}]$$
allows us to prove that $k_{-1}y_{-2n} = x_{-4n+1}\alpha_{2n-2}$ and that $k_{-1}\beta_{4n-2}$ and $k_{-1}\alpha_{2n-2}$
don't vanish.

\rightline{$\square$}

\vspace{3mm}

Let us now explain how to compute the non zero differentials on $\mathbb{E}^1_{*,*}(\mathcal{MS})(\pi)$.
The module $\mathbb{E}^{\infty}_{-4n+1,2n-2}[ev(0)]$ is $\mathbb{Z}/2 \mathbb{Z}$. Then, we must have $$d_1(kiT)= 2ac.$$ Since $d_1(iT)=0$ by comparison with the Cohen-Jones-Yan spectral sequence, we must have $$d_1(kT)=2c.$$ 
The multiplication by $c \in \mathbb{E}^1_{0, 2n-2}(\mathcal{MS})(\pi)$ gives the differentials on $c^{\kappa}kT=l^{\kappa}kT$  and on  $c^{\kappa}kiT=l^{\kappa}kT$, $\kappa \in \mathbb{N}$.
The comparison with the Cohen-Jones-Yan spectral sequence allows us to conclude that they are the only non zero differentials from the first column of $\mathbb{E}^*_{*,*}(\mathcal{MS})(\pi) $. 
We deduce the other differentials by multiplying by $T$.

\begin{center}
\includegraphics[width=11cm]{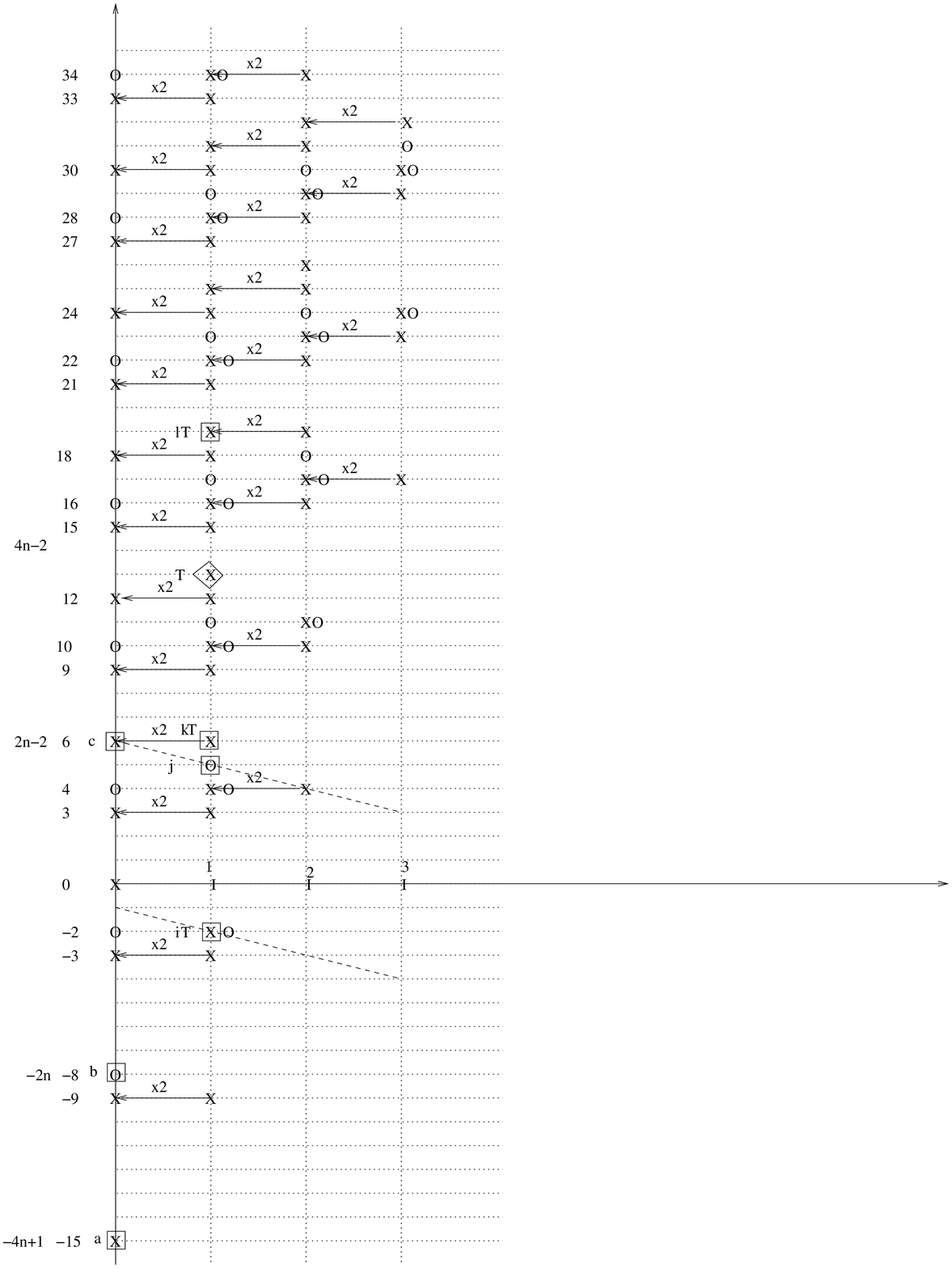} 
\end{center}

The new generators are on the two first columns of $\mathbb{E}^2_{*,*}(\mathcal{MS})(\pi)$. For degree reasons there could not have non zero differentials on this page and on the others. Thus we have  $\mathbb{E}^2_{*,*}(\mathcal{MS})(\pi) = \mathbb{E}^{\infty}_{*,*}(\mathcal{MS})(\pi)$.

\begin{center}
\includegraphics[width=13cm]{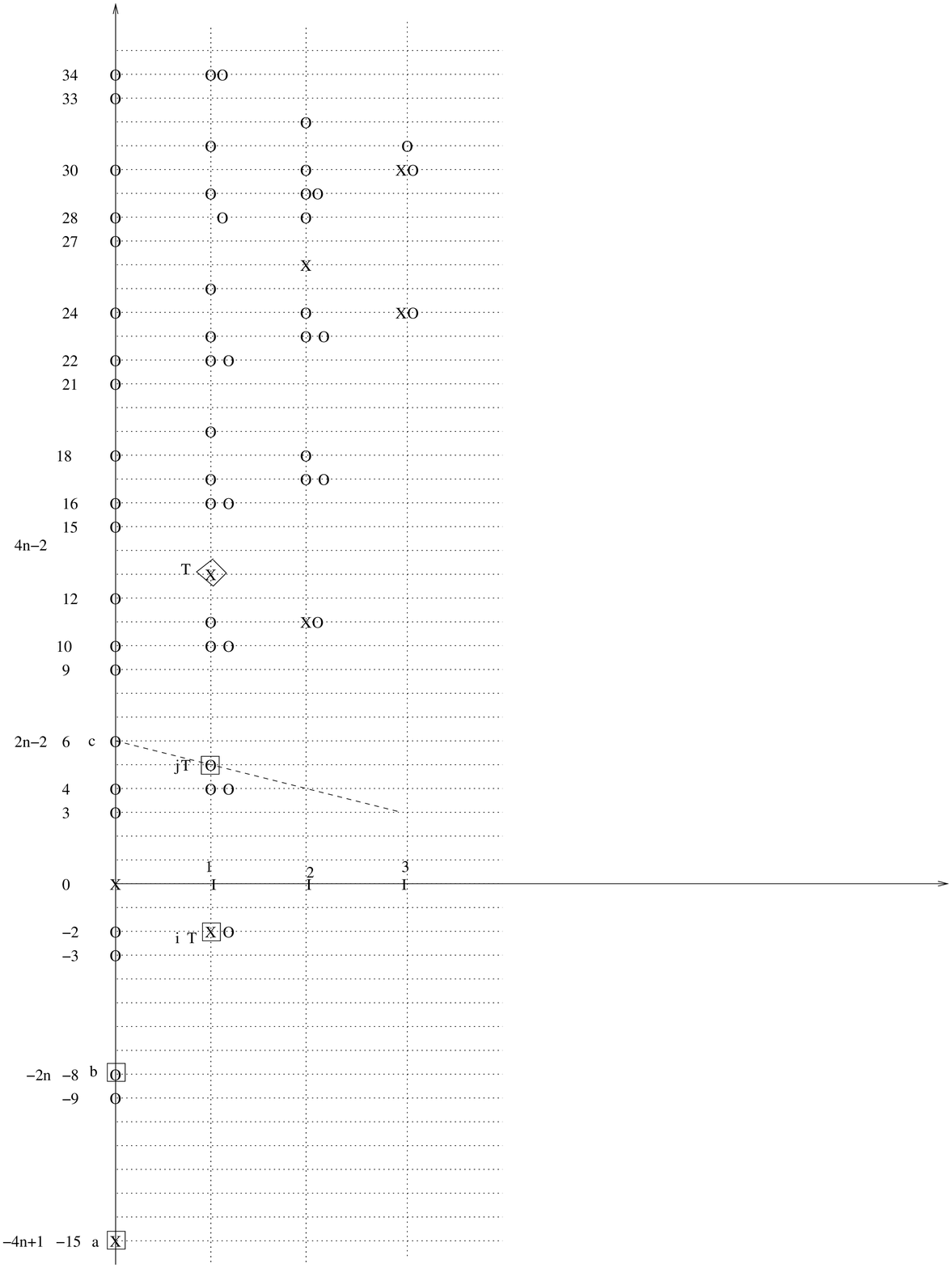} 
\end{center}

\subsection{Solve the linear extension issues.}
Now that we know the $E^{\infty}$-terms of the two spectral sequences, we have to compare them to solve their extension issues.
In the Cohen-Jones-Yan spectral sequence we have two linear extensions issues. They are located in total degree $-1$ and in degree $2n-2$.

\begin{proposition} \label{linext}
The module $\mathbb{H}_{2n-2}(\mathcal LUS^{2n})$ is equal to $\mathbb{Z}_2 \oplus \mathbb{Z}_2$.
\end{proposition}
\noindent{Proof.}
To solve this extension issue, we use the mulplicative structure of the spectral sequences. By our preceding computations we have to decide if $\mathbb{H}_{2n-2}(\mathcal LUS^{2n}) \simeq \mathbb{Z}_2 \oplus \mathbb{Z}_2$ or $\mathbb{Z}_4$.
Let $x_{2n-2}$ be an element of  $\mathbb{H}_{2n-2}(\mathcal LUS^{2n})$ that represents the element $\alpha_{2n-2}$ in the Cohen-Jones-Yan spectral sequence and let $x_{-4n+1}$ be the element of lowest degree in $\mathbb{H}_{*}(\mathcal LUS^{2n})$. Then, the product $x_{2n-2}  x_{-4n+1}$ lies in $\mathbb{E}^{\infty}_{0,-2n-1}(\mathcal{MS})(\pi)$. This proves that $x_{2n-2}$ lies in $F_0\mathbb{H}_{2n-2}(\mathcal LUS^{2n})$ according to the Morse-Serre filtration and that this generator can not be of $4$-torsion but of $2$-torsion.
\rightline{$\square$}

\subsection{Solve the extension issues of algebra.}

\begin{proposition} \label{extalg}
There is no extension issues of algebra on the Cohen-Jones-Yan spectral sequence $\mathbb{E}^{\infty}_{*,*}[ev(0)]$.
\end{proposition}
\noindent{Proof.}

Now that we have solve the linear extension issues, we can write that $x_{-4n+1}$, $y_{-2n}$, $\alpha_{2n-2}$, $\beta_{4n-2}$ and $k_{-1}$ as elements of $\mathbb{H}_*(\mathcal LUS^{2n})$.

\begin{enumerate}
\item From the Morse-Serre filtration, we know that $x_{-4n+1}$, $y_{-2n}$ and $\alpha_{2n-2}$ lie in $\mathbb{E}^{\infty}_{0,*}(\mathcal{MS})(\pi)$. Then, there is no extension issue of algebra on products between them.

\item $x_{-4n+1}k_{-1} = 0$ for degree reasons.

\item There is no extension issues concerning the products $x_{-4n+1}\beta_{4n-2}$, $y_{-2n}k_{-1}$, $y_{-2n}\beta_{4n-2}$ and $k_{-1}^2$ because of there position in the filtration.

\item On the Cohen-Jones-Yan spectral sequence, we see that there is a free module in the filtration preceding the product $\beta_{4n-2}k_{-1}$. Since there is no $2$-torsion element in this module, there is no ambiguity.

\item The possible ambiguity concerning $\beta_{4n-2}^2$ is solved by computing $\mathbb H_*(\mathcal LUS^{2n},\mathbb Q)$. We use the fact that over the rational numbers $US^{2n}$ is homotopy equivalent to $S^{4n-1}$ (the loop product is a rational homotopy invariant).

\item The two last extension issues of algebra namely $\alpha_{2n-2}k_{-1}$ and  $\alpha_{2n-2}\beta_{4n-2}$ are solved by using the Morse-Serre spectral sequence where there are no ambiguities for these products.

\end{enumerate}

\rightline{$\square$}

\section{Computation of $\mathbb{H}_*(\mathcal LUS^{2n+1})$, $n \geq 2$.}
\noindent{\bf Strategy of the proof.}
We compute $\mathbb{E}^{\infty}_{*,*}(\pi)$,  the $E^{\infty}$ term of the $(2n+1, 2n)$-regraded Serre spectral sequence associated to the fibration $\pi: \mathcal LUS^{2n+1} \to \mathcal LS^{2n+1}$ in 5.1. There is no linear extension issue. Then, by computing the Cohen-Jones-Yan spectral sequence $\mathbb{E}^*_{*,*}(ev(0))$ in 5.2, and the Morse Serre spectral sequence $\mathbb{E}^*_{*,*}(\mathcal{MS}(\pi))$ in 5.3, we solve the extension issues of algebra from $\mathbb{E}^{\infty}_{*,*}(\pi)$ in 5.4.

\subsection{The regraded Serre spectral sequence $\mathbb{E}^{*}_{*,*}(\pi)$.}
We compute the last page of the $(2n+1,2n)$-regraded Serre spectral sequence of the fibration $\pi:\mathcal LUS^{2n+1} \to \mathcal LS^{2n+1}$. Then, we deduce the graded module structure of $H_*(\mathcal LUS^{2n+1})$.

\begin{proposition} \label{m}
The last page of $\mathbb{E}^{*}_{*,*}(\pi)$ namely $\mathbb{E}^{\infty}_{*,*}(\pi)$ is isomorphic to $$\mathbb{H}_*(\mathcal LS^{2n+1}) \otimes \mathbb{H}_*(\mathcal LS^{2n}) $$ that is isomorphic to
$$\mathbb{Z}[x_{-2n-1}, v_{2n}]/(x_{-2n-1}^2) \otimes \mathbb{Z}[y_{-2n}, u_{4n-2}, \theta_{-1}]/(y_{-2n}^2, y_{-2n}\theta_{-1}, \theta_{-1}^2, 2y_{-2n}u_{4n-2}).$$

\end{proposition}

Proof.

Let us determine the $(2n+1, 2n)$-regraded Serre spectral sequence of the fibration $\pi: \mathcal LUS^{2n+1} \to \mathcal LS^{2n+1}$.
The Chas and Sullivan loop product povides this spectral sequence with a multiplicative structure (see \cite{L}).
Since $\mathbb{H}_*(\mathcal LS^{2n+1})$ has no torsion, $\mathbb{E}^{2}_{*,*}(\pi) \cong \mathbb{H}_*(\mathcal LS^{2n+1}) \otimes \mathbb{H}_*(\mathcal LS^{2n})$.
There are five generators of algebra : $x$ of degree $-(2n+1)$, $y$ of degree $-2n$,  $u$ of degree $4n-2$, $v$ of degree $2n$ and $\theta$ of degree $-1$.
This generators are represented by squares on the following diagram.

\begin{center}
\includegraphics[width=16cm]{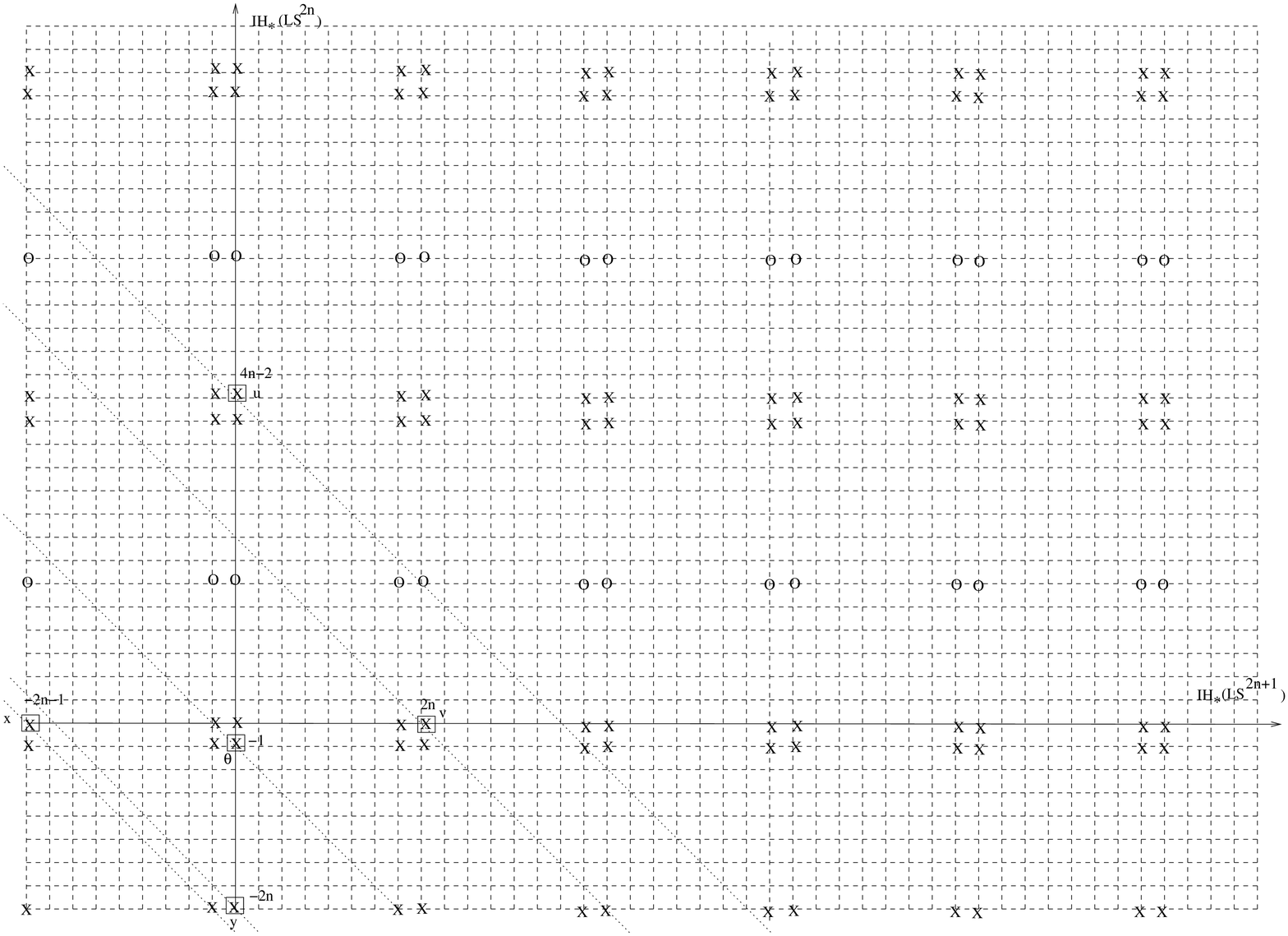} 
\end{center}

For degree reasons the only non zero differential starting from a generator is $d_{2n}$. More precisely, we only have to compute $d_{2n}(y)$.
Since there is a section $US^{2n+1} \to \mathcal LUS^{2n+1}$ and since that $\mathbb{H}_*(US^{2n+1}) = \mathbb{H}_*(S^{2n+1}) \otimes \mathbb{H}_*(S^{2n})$, we have that $d_{2n}(y)=0$.
This proves that this spectral sequence collapse at the $E^2$-level and that the following figure represents $\mathbb{E}^{\infty}_{*,*}(\pi)$.

\rightline{$\square$}

By an easy check we have the following proposition.

\begin{proposition} \label{w}
There is no linear extension issue in $\mathbb{E}^{\infty}_{*,*}(\pi)$. Then, $H_*(LUS^{2n+1})$ is isomorphic to $H_*(\mathcal LS^{2n+1}) \otimes H_*(\mathcal LS^{2n})$ as a graded module.
\end{proposition}

\subsection{The Cohen-Jones-Yan spectral sequence, $\mathbb{E}^*_{*,*}(ev(0))$.}
Since we know the module structure of $H_*(\mathcal LUS^{2n+1})$, we can deduce all the differentials of the Serre spectral sequence $\mathbb E_{*,*}^{*}(ev(0))$ of the fibration $ev(0): \mathcal LUS^{2n+1} \to US^{2n+1}$.

\begin{proposition} \label{x}
The algebra $\mathbb{E}^{\infty}_{*,*}(ev(0))$ is isomorphic to  $$\mathbb{Z}[a_{-2n-1}, b_{-2n}, c_{-1}, d_{2n}, e_{4n-2}]/(a_{-2n-1}^2, b_{-2n}^2, c_{-1}^2, c_{-1}b_{-2n}, 2e_{4n-2}b_{-2n}).$$
\end{proposition}

Proof.
First recall that $H_*(\Omega US^{2n+1}) \cong H_*(\Omega S^{2n+1}) \otimes H_*(\Omega S^{2n})$.
Since $H_*(\Omega US^{2n+1})$ is without torsion, then 
$$\mathbb{E}^2_{*,*}(ev(0)) \cong \mathbb{H}_*(US^{2n+1}) \otimes H_*(\Omega US^{2n+1})$$
$$\cong \mathbb{Z}[a_{-2n-1}, b_{-2n}, d_{2n}, f_{2n-1}]/(a_{-2n-1}^2, b_{-2n}^2).$$
These generators are represented by squares on the following diagram.

\begin{center}
\includegraphics[width=14cm]{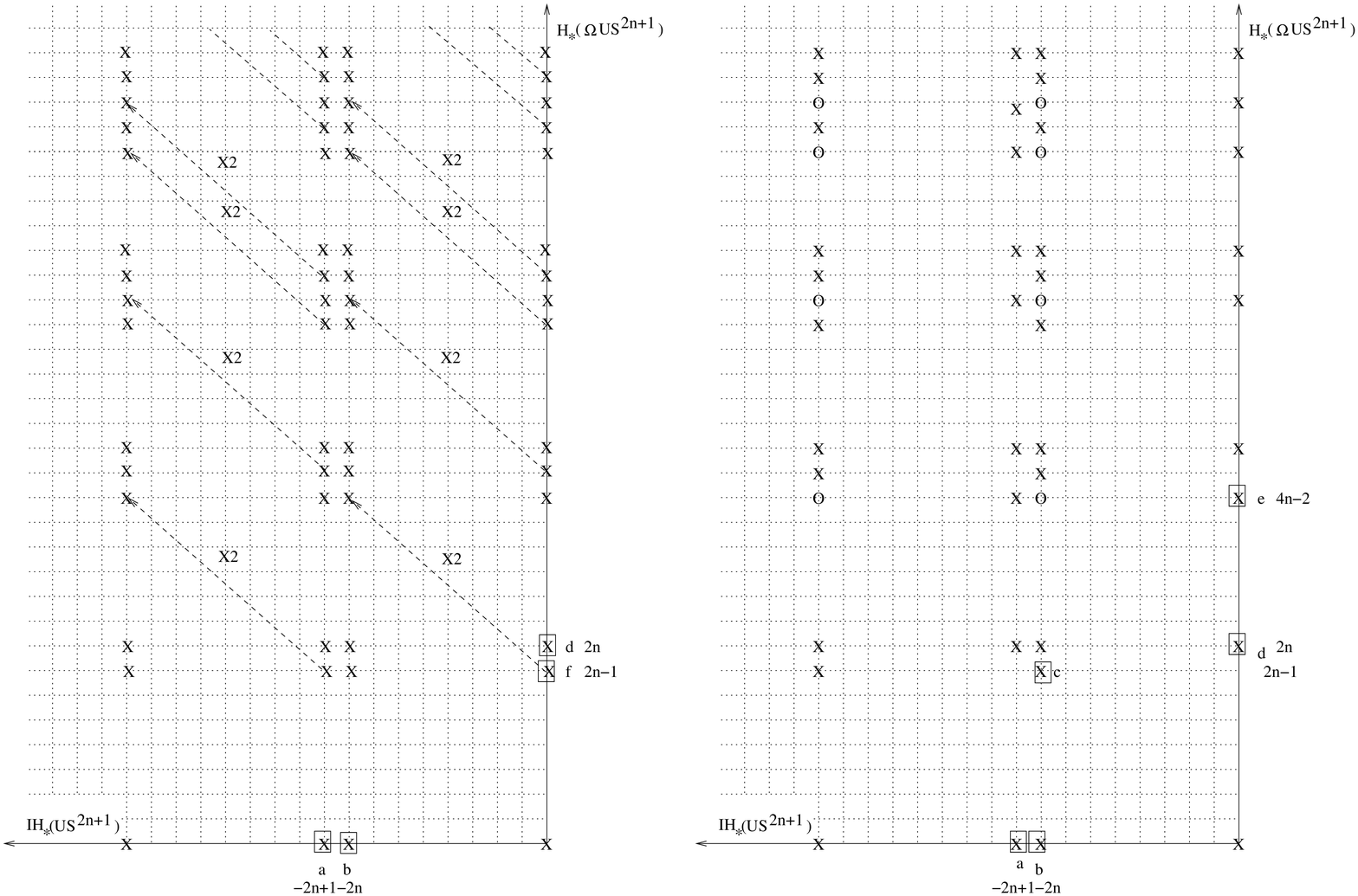} 
\end{center}

Since we know the linear structure of $\mathbb{H}_*(LUS^{2n+1})$, we can deduce all the differentials and compute $\mathbb{E}^{\infty}_{*,*}(ev(0))$ together with its structure of algebra.

\rightline{$\square$}

\subsection{The Morse Serre spectral sequence, $\mathbb{E}^{*}_{*,*}(\mathcal{MS})(\pi)$.}
As for the even case, we determine $\mathbb{E}^{1}_{1,*}(\mathcal{MS})(\pi)$ by computing the $(4n+1, 2n)$-regraded Serre spectral sequence of the fibration $b: \widetilde{US^{2n+1}} \to US^{2n+1}$.

\begin{proposition} \label{p}
The subalgebra $\mathbb{E}^1_{p,*}(\mathcal{MS})(\pi)$, $p \geq 1$, is isomorphic to
$$\mathbb{Z}[\tau, \upsilon, \phi, \chi, \psi]/(\tau^2, \upsilon^2, \chi^2, \phi^2, 2 \psi \phi)[T]_{\geq 1}$$
with $deg(\tau)=-2n-1$, $deg(\upsilon)=-2n$, $deg(\phi)=-2n$, $deg(\chi)=-1$, $deg(\psi)=4n-2$ and $bideg(T)= (1,4n-1)$.

\end{proposition}

\noindent{Proof.}
We compute the $E^2$ term of the $(4n+1, 2n)$-shifted Serre spectral sequence associated to the fibration $b: \widetilde{US^{2n+1}} \to US^{2n+1}$.

\begin{center}
\includegraphics[width=7cm]{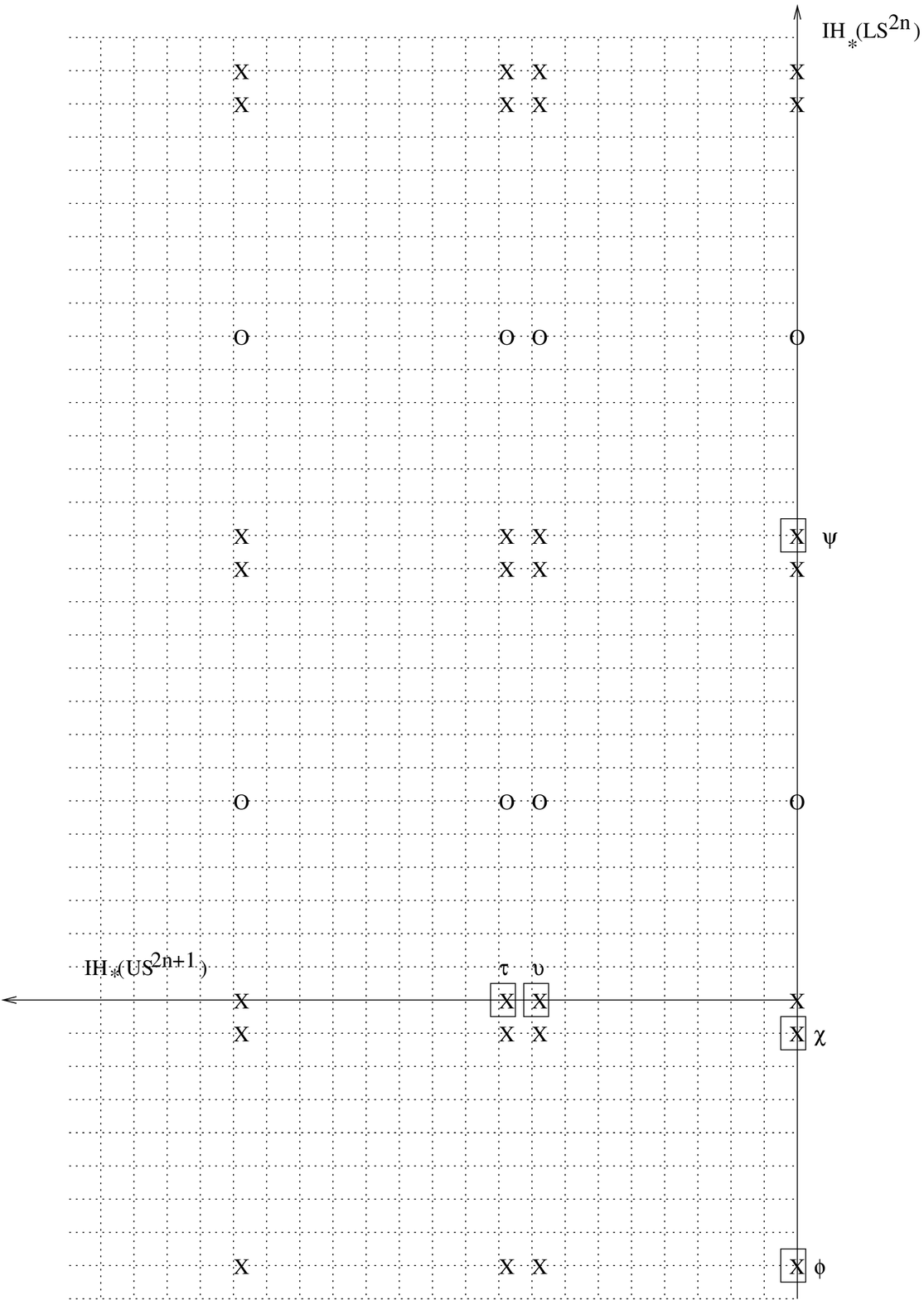} 
\end{center}

An easy inspection and the use of the section of $\pi: \widetilde{\Sigma_1} \to \Sigma_1$ as for the even dimensional case (section 4.1), shows that all the differentials of this spectral sequence must vanish.
\rightline{$\square$}

We deduce that $\mathbb{E}^{*}_{*,*}(\mathcal{MS})(\pi)$ collapses at the $E^1$ level and we get $\mathbb{E}^{\infty}_{*,*}(\mathcal{MS})(\pi)$.

\begin{center}
\includegraphics[width=3cm]{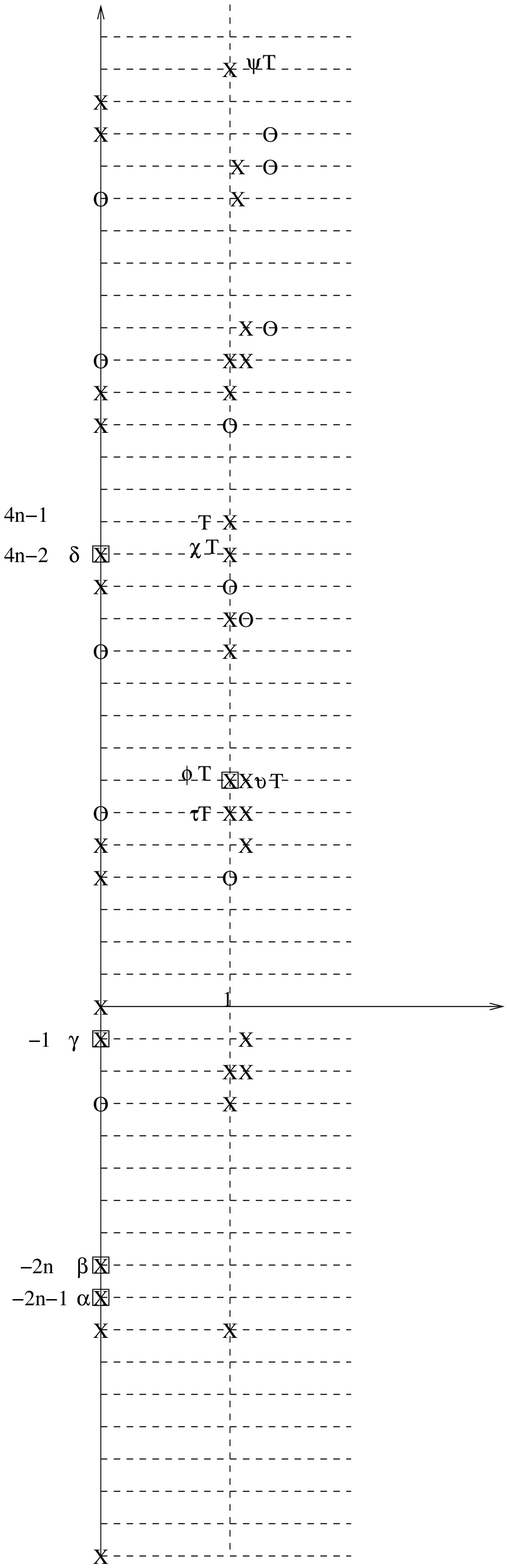} 
\end{center}

\begin{proposition}
The algebras $\mathbb{H}_*(\mathcal LUS^{2n+1})$ and  $\mathbb{E}^{\infty}_{*,*}(\pi)$ are isomorphic. 
\end{proposition}

\noindent{Proof.}
There may have ambiguities on the products of generators on $\mathbb{E}^{\infty}_{*,*}(\pi)$ for the following products.
\begin{enumerate}
 \item The products involving $\theta_{-1}$, $y_{-2n}$ and $u_{4n-2}$. These generators lie in $\mathbb{E}^{\infty}_{0,*}(\mathcal{MS})(\pi)$ so there is no multiplicative extension issue.

\item The product $v_{2n}\theta_{-1}$. We can drop the extension issue by considering the Cohen-Jones-Yan spectral sequence $\mathbb{E}^*_{*,*}(ev(0))$.

\item The last extension issue concerning the product $v_{2n}\theta_{-1}$ can be solved by considering the multiplicative structure of the Morse-Serre spectral sequence and to be more precise the $\mathbb{E}^{\infty}_{0,*}(\pi)$-module structure on $\mathbb{E}^{\infty}_{1,*}(\pi)$.
\end{enumerate}

\rightline{$\square$}

\end{document}